
\magnification=\magstep1
\input amstex
\documentstyle{amsppt}
\hsize=5in
\vsize=7.8in
\strut
\vskip5truemm
\font\small=cmr8
\font\itsmall=cmti8

\def\smallarea#1{\par\begingroup\baselineskip=10pt#1\endgroup\
par}
\def\abstract#1{\begingroup\leftskip=5mm \rightskip 5mm
\baselineskip=10pt\par\small#1\par\endgroup}

\def\refer#1#2{\par\begingroup\baselineskip=11pt\noindent\left
skip=8
.25mm\rightskip=0mm
 \strut\llap{#1\kern 1em}{#2\hfill}\par\endgroup}

\nopagenumbers
\centerline{\bf JUMPS OF THE ETA-INVARIANT}
\vskip 25.3pt

\centerline{MICHAEL S. FARBER
\footnotemark"$^1$"}
\footnotetext"$^{1,2}$"{The research was supported
by grant No.~88-00114 from the United States - Israel
Binational
Science
Foundation (BSF), Jerusalem, Israel}
\footnotetext"$^{2}$"{Research supported in part
by National Science Foundation}
\smallarea{%
\centerline {\itsmall School of Mathematical Sciences,}
\centerline {\itsmall Tel-Aviv University,}
\centerline {\itsmall Tel-Aviv 69978, Israel}
}
\medskip
\centerline {JEROME P. LEVINE\footnotemark"$^2$"}
\smallarea{%
\centerline{\itsmall Department of Mathematics,}
\centerline{\itsmall Brandeis University,}
\centerline{\itsmall Waltham, 02254 Mass USA}}

\vskip 0.5in
fibred
knot,


\input amstex
\define\C{{\Bbb C}}
\define\R{{\Bbb R}}
\define\RR{{\Cal R}}
\define\Z{{\Bbb Z}}
\define\Q{{\Bbb Q}}

\define\OO{{\Cal O}}
\define\UU{{\frak U}}
\define\T{{\Cal T}}

\define\p{{\frak p}}
\define\Res{\operatorname{Res}}
\define\Hom{\operatorname{Hom}}
\define\ch{\operatorname{ch}}
\define\ev{\operatorname{ev}}
\define\index{\operatorname{index}}

\define\End{\operatorname{End}}
\define\Har{\operatorname{Har}}
\define\Diff{\operatorname{Diff}}
\define\M{{\Cal M}}
\define\A{{\Cal A}}

\define\BB{{\frak B}}
\define\VV{{\Cal V}}
\define\m{{\frak m}}
\define\F{{\Cal F}}
\define\FF{{\frak F}}
\define\sign{\operatorname{sign}}
\define\tr{\operatorname{tr}}
\define\im{\operatorname{im}}
\redefine\H{\Cal H}

\define\E{\Cal E}
\redefine\L{{\Cal L}}
\def\<{\langle}
\def\>{\rangle}
\documentstyle{amsppt}

\define\ABP{1}
\define\APS{2}
\define\BiC{3}
\define\BGV{4}
\define\Bla{5}
\define\BC{6}
\define\Bo{7}
\define\CS{8}
\define\Dai{9}
\define\Fa{10}
\define\FKK{11}
\define\Gan{12}
\define\Gilkey{13}
\define\Godement{14}
\define\Kato{15}
\define\Kirk{16}
\define\Kob{17}
\define\Levin{18}
\define\Levine{19}
\define\Levinee{20}
\define\Mathai{21}
\define\Ne{22}
\define\No{23}
\define\Pal{24}
\define\Rud{25}
\define\Wa{26}
\define\War{27}
\define\We{28}

\topmatter
\endtopmatter
\TagsOnRight
\TagsAsMath
\document

In \cite{\APS}, Atiyah, Patodi and Singer introduced an
invariant
$\eta_D$ of any self-adjoint elliptic differential operator
$D$ on
an
odd-dimensional oriented closed manifold $M$, in order to
prove an
index theorem for manifolds with boundary. For the germinal
case of
the
\lq\lq signature operator\rq\rq the relevant $D$ is $\pm(\ast
d-
d\ast)$,
where the Hodge duality operator $\ast$ is determined by the
Reimannian
metric on $M$. They consider, more generally, the signature
twisted
by a flat
connection $\nabla$ on a Hermitian vector bundle $\E$, or
equivalently,
a unitary representation $\alpha$ of $\pi_1(M)$. Then $D$ is
replaced by
$\pm(\ast \nabla-\nabla\ast)$ and the corresponding eta-
invariant is
denoted $\eta_\nabla$. An important observation is that the
\lq\lq reduced\rq\rq invariant $\eta_\nabla-k\eta_d$ (where
$k=\dim
\E$)
is a \lq\lq topological\rq\rq (more precisely, a $C^\infty$)
invariant
of $(M,\alpha)$ - the $\rho$-invariant.

The eta-invariant $\eta_\nabla$, considered as a real valued
function of a
flat connection $\nabla$, demonstrates two different
phenomena.
Firstly, it has {\it integral jumps\/}, known also as the
{\it spectral flow\/}, playing a central role in modern low-
dimensional
topology. And, secondly, it varies {\it smoothly} if
$\eta_\nabla$
is
considered modulo integers, i.e. as a function (called {\it
the
reduced eta-invariant}) with values in $\R/\Z$.

The results of this paper contribute to the understanding of
these phenomena. Our
main result, Theorem 1.5, states that the jumps (or,
equivalently,
the infinitesimal spectral flow) can be calculated {\it
homologically\/}
by means of a linking form, constructed directly in terms of
{\it deformations of the monodromy
representations of the fundamental group}. More specifically,
we show that given an analytic path of unitary
representations, one
may
explicitly construct a local coefficient system over the
manifold
and
an {\it algebraic} linking form on the homology of this local
system
such that the signature
invariants of the linking form determine the spectral flow
completely.

To prove the above mentioned result we study first a more
general
situation
of analytic families $D_t$ of arbitrary elliptic self-adjoint
operators
on closed manifolds. We show that any such family can be
viewed as a
single operator $\tilde D$ acting on the space of germs of
analytic
curves of
smooth sections; such operator $\tilde D$ determines, in a
canonical
way,
a linking pairing whose signatures measure the jumps of the
eta-
invariant
of the family. This linking form we call {\it analytic} in
order to
distinguish it from the algebraic linking form mentioned
above.
In order to prove our principal result, Theorem 1.5, we need
to find a relation between
these two linking forms (they are not isomorphic, although
have
the same signatures); this is done
with the use of parametrized Hodge decomposition, described in
\S 4,
and
a version of the De Rham theorem for the germ-complex, cf.
\S5.
These are
the main ingredients of the proof.

A recent preprint \cite{\Kirk} of Paul Kirk and Eric Klassen
also
addresses
the problem of homological computation of the spectral flow.
They
proved
that the contribution of the "first order terms" to the
spectral
flow
is equal to the signature of a quadratic form given as a cup-
product
(compare Corollary 3.15).
Our Theorem 1.5 gives a more precise answer; it describes the
contributions of
terms of all orders.

Theorem 1.5 has some ideological similarity with the results
of
X.Dai \cite{\Dai},
who studied adiabatic limit of the eta-invariant in the case
when
the Dirac
operators along the fibers have nontrivial kernels, forming a
vector
bundle.
Before him the adiabatic limit formula was obtained by
J.-M.Bismut and J.Cheeger \cite{\BiC} under the assumption
that the
kernels
along the fibers are trivial.
X.Dai proved in \cite{\Dai}
that the adiabatic limit formula for the
signature operator contains an additional topological
invariant
$\tau$ which has a twofold characterization:
it measures nonmultiplicativity of the signature and is equal
to the
sum of
the signatures of pairings determined on the terms of the
Leray
spectral
sequence.

We apply our theorem about jumps of the eta-invariant to study
the problem of homotopy invariance of the
Atiyah--Patodi--Singer $\rho$-invariant.
We give an intrinsic
homotopy theoretic definition of the $\rho$-invariant, up to
indeterminacy
in the form of a locally constant function on the space of
unitary
representations, which is
zero at the trivial representation. The proof relies on some
auxiliary
results (mostly known) about variation of the mod $\Z$-
reduction of
the
eta-invariant. We conjecture that this indeterminacy is
rational-
valued.

A consequence of our results is that the $\rho$-invariants of
homotopy
equivalent manifolds differ by a rational-valued, locally
constant
function
on every component of the representation space containing some
representation
which factors through a group satisfying the Novikov
conjecture. But
Shmuel Weinberger shows that such representations are actually
dense
in the
representation space and so the difference is rational valued
everywhere.
We are grateful to Weinbergerfor including his proof as an
Appendix
to
our paper.

We are also grateful to M.Braverman, J.Hillman, P.Kirk and
S.Shnider
for
stimulating discussions. Many thanks to V.Matsaev for his help
with proof of Lemma 5.2 and for numerous consultations
concerning
analytic subtleties.

\heading 1. Deformations of the monodromy representations and
jumps
of the eta-invariant\endheading

1.1. Let $M$ be a compact oriented Riemannian manifold of
odd dimension $2l-1$
and $\E$ a {\it flat Hermitian vector
bundle} of rank $m$ over $M$. This means that
(1) a Hermitian metric has been specified on
each fibre $\E_x$ which
varies smoothly with $x \in M$; (2) there is given a covariant
derivative
$$\nabla : A^k(M;\E)\to A^{k+1}(M;\E), \ k=0,1, \dots\tag1$$
acting on the space of $C^\infty$-forms on $M$ with values in
$\E$;
(3) $\nabla$ is {\it flat}, i.e.
$\nabla ^2=0;$
and (4) the covariant derivative $\nabla$ is
{\it compatible} with the Hermitian
structure on $\E$; the latter can also be expressed by saying
that
the
Hermitian metric on $\E$ is {\it flat}.

In this situation Atiyah, Patodi and Singer \cite{\APS} have
defined the following first order differential operator
acting on forms of even degree $\phi \in A^{2p}(M,\E)$ by
$$B\phi=i^l(-1)^{p+1}(\ast \nabla -\nabla \ast)\phi,
\ \ \ B: A^{ev}(M;\E) \to A^{ev}(M;\E)\ \tag2$$
where the star denotes the Hodge duality operator.
The operator $B$ is elliptic and self-adjoint. To any such
operator
Atiyah, Patodi and Singer in \cite{\APS} assigned a numerical
invariant,
$\eta(B)$, called the eta-invariant which plays a crucial role
in
the index
theorem for manifolds with boundary.
Recall that the eta-invariant $\eta(B)$ is defined as follows.
Consider the eta-function of $B$:
$$\eta_B(s)=\sum_{\lambda\ne 0}\sign(\lambda)|\lambda|^{-s},$$
where $\lambda$ runs over all eigenvalues of $B$. It follows
from
the
general theory of elliptic operators
that for large $\Re (s)$ this formula defines a holomorphic
function
of $s$
which has a meromorphic continuation to the whole complex
plane.
Atiyah, Patodi and Singer proved in \cite{\APS} that the eta-
function is
holomorphic at $0\in \C$ (in the general situation this was
proven
by
P. Gilkey \cite{\Gilkey}, Theorem 4.3.8).
The eta-invariant of $B$ is then defined
as the value
of the eta-function at the origin $\eta_B(0)$.

1.2. Suppose now that the covariant
derivative $\nabla$ is being {\it deformed}. By this
we understand that there is given an
{\it analytic family} of differential operators
$$\nabla_t: A^k(M;\E) \to A^{k+1}(M;\E), \ \ k=0, 1,
\dots\tag3$$
where the parameter $t$ varies in an interval around zero
$(-\epsilon, \epsilon)$
such that:

(i) for every value of $t$ the operator $\nabla_t$ is a
covariant derivative on the vector
bundle $\E$ having
curvature $0$ (i.e. $\nabla_t^2=0$) and the Hermitian metric
on $\E$
is flat with respect to every $\nabla_t$;

(ii) for $t=0$ the operator $\nabla_0$ coincides with the
original
covariant derivative $\nabla$.

The analyticity of the family of connections $\nabla_t$ we
understand as
follows. Represent $\nabla_t=\nabla +\Omega_t$ where
$\Omega_t\in A^1(M;\End(\E))$; then the curve $t\mapsto
\Omega_t$ is
supposed to be analytic with respect to any Sobolev norm; cf.
section 3.3
below.

The corresponding self-adjoint operators $B_t$, constructed
using
the connections $\nabla_t$
as explained above, will also be functions
of the parameter $t$. Let $\eta_t$ denote the eta-invariant of
$B_t$.
It is shown in
\cite{\APS} that $\eta_t$ may have only simple discontinuities
--
integral jumps occuring when some eigenvalues of $B_t$ cross
zero.
In other words,
the limits
$$\lim_{t\to +0}\eta_t=\eta_+,\ \ \ \lim_{t\to
-0}\eta_t=\eta_-\tag4
$$
exist and the discontinuities (jumps)
$\eta_+ -\eta_0$ and $\eta_- -\eta_0$ are integers; here
$\eta_0$
denotes
the eta-invariant of the operator $B_0$. Atiyah, Patodi and
Singer \cite{\APS} show that the integer
$\eta_+ -\eta_-$ has the meaning of {\it infinitesimal
spectral
flow};
the integers
$\eta_+ -\eta_0$ and $\eta_- -\eta_0$ have the meaning of
{\it infinitesimal half flows}.

Our aim in the present paper is to compute these integral
jumps
$\eta_+ -\eta_0$ and $\eta_- -\eta_0$ in terms of some {\it
homological\/}
invariants constructed by means of the {\it germ of the
deformation}
$\nabla_t$. More precisely, we will express these jumps in
terms of
the germ of the {\it deformation of the monodromy
representation\/}.

1.3. Fix a base point $x\in M$. Given a flat
Hermitian connection $\nabla$
on $\E$ there is the corresponding monodromy representation
$$\rho: \pi=\pi_1(M,x)\to  U(\E_x)$$
where $\E_x$ is the fibre above $x \in M$ and $U(\E_x)$
denotes
the unitary group. We can express
this by saying that $\E_x$ has the structure of a left
$\C[\pi]$-
module,
where for $g\in \pi$ and $v\in \E_x$ the product $g\cdot
v\in\E_x$
is the
value of the flat section obtained from $v\in \E_x$ by
parallel displacement along
a loop representing $g$.

When the flat connection is being deformed,
the corresponding monodromy
representation is deformed as well. Thus, we obtain from the
family
$\nabla_t$ an {\it analytic one-parameter family
of left $\C[\pi]$-module structures on} $\E_x$. In other
words, we
obtain
a family of maps
$$\mu_t: \pi\times\E_x\to \E_x, \ \ (g,v)\mapsto g\cdot_tv, \
\
g\in\pi,v\in\E_x$$
where $t\in (-\epsilon,\epsilon)$, such that for any $t$ the
map
$(g,v)\mapsto g\cdot_tv$
is a linear unitary action of $\pi$ on $\E_x$, depending
analytically on
$t$. (The analyticity follows from Lemma 5.2 below).

Let $\OO$ denote the ring of germs of
complex valued holomorphic
functions
$f:\ (-\epsilon,\epsilon)\to \C$
at the origin.
An element of $\OO$ can
also be represented by a power series
$$f(t)=\sum_{n\ge0}a_nt^n, \ \ a_n\in\C$$
having a non-zero radius of convergence.
We will consider $\OO$ together
with the involution which is induced by complex conjugation on
$\C$;
it has the property that $\overline t =t$ (i.e. $t$ is real).

Let $\OO\E_x$ be the set of germs of holomorphic curves in
$\E_x$:\
\
$\alpha:(-\epsilon,\epsilon)\to \E_x.$ It is a left $\OO$-
module
where
$$(f\cdot\alpha)(t)=f(t)\cdot\alpha(t), \ \ f\in\OO,\quad
\alpha\in\OO\E_x.$$
It is clear that $\OO\E_x$ is free of rank $m=\dim\E_x$ over
$\OO$.
The
Hermitian metric on the fibre $\E_x$ defines (by pointwise
multiplication)
the following bilinear map
$$\<\ ,\ \>:\OO\E_x\times \OO\E_x\to \OO\tag5$$
which is Hermitian, non-degenerate and $\OO$-linear
with respect to the first variable.

Given a deformation of the monodromy representation as above
consider
the following map
$$\pi\times\OO\E_x\to\OO\E_x,\ \ (g,\alpha)\mapsto
g\cdot\alpha$$
where $g\in\pi,\ \alpha\in \OO\E_x$ and $g\cdot\alpha$ denotes
the germ of the following curve
$$t\mapsto (g\cdot\alpha)(t)=g\cdot_t(\alpha(t))\ \ \in \E_x$$
This map defines a left $\OO[\pi]$-module structure on
$\OO\E_x$.

We arrive at the conclusion that
{\it \lq\lq an analytic one parameter family of left
$\C[\pi]$-module structures on $\E_x$\rq\rq} can be understood
as a
left module
$\VV=\OO\E_x$ over the ring $\OO[\pi]$ having the following
properties:

(a) $\VV$ is free of rank $m=\dim \E_x$ over $\OO$ and is
supplied with a Hermitian form
$$\<\ ,\ \>:\VV\times \VV\to \OO;$$

(b) the form $\<\ ,\ \>$ is non-singular and $\OO$-linear
in the first variable and
anti-linear in the second variable;

(c) $\<\lambda v,w\>=\<v,\overline\lambda w\>$ for $\lambda
\in
\OO[\pi],
\ v, \ w \in \VV$. Here the involution acts on $g\in \pi$ as
$\overline g=g^{-1}$;

(d) If $\m$ denotes the maximal ideal of $\OO$ then
there is an isomorphism
$$\VV/\m \VV \thickapprox \E_x$$ of $\C[\pi]$-modules
such that the form $\<\ ,\ \>$ reduces under this isomorphism
to
the original Hermitian form on the fibre $\E_x$.

The pair consisting of the $\OO[\pi]$-module $\VV$ and the
form $\<\
,\ \>$
on it will be referred to as
{\it deformation of the monodromy representation\/}. Thus, a
deformation
$\nabla_t$ of a flat Hermitian bundle determines a deformation
of
the
monodromy representation. In 5.3 we will describe $\VV$ as a
monodromy of a
locally flat preseaf over $M$.

1.4. Now we will show that some standard homological
constructions
(using
Poincare duality) lead to certain
{\it linking forms\/} constructed by means of
the deformation of the monodromy representation. The linking
form
described
here we will sometimes call {\it homological\/} or {\it
algebraic\/}
in
order to distinguish it from another linking form also
determined by
the deformation which we will call {\it analytic\/}; it is
constructed in
section 3.

Consider the cohomology of the manifold $M$
with local coefficient system defined by
$\VV$; this we understand as
$$H^{\ast}(\Hom
_{\Z[\pi]}
(C_{\ast}(\tilde M),\VV))$$
and will denote it by $H^{\ast}(M;\VV).$
Here $\tilde M$ is the universal cover of $M$ and the group
$\pi$
acts
on $\tilde M$ from the left by covering translations.
Note that $H^{\ast}(M;\VV)$ is a finitely generated $\OO$-
module.
Since
$\OO$ is a principal ideal domain the module $H^k(M;\VV)$ can
be
represented
as the sum
$$H^k(M;\VV)=F^k\oplus T^k$$
of its $\OO$-torsion $T^k=T^k(M)$ and the free part
$F^k=H^k(M;\VV)/T^k$
for every $k=0, 1, \dots$. Note that $T^k(M)$ is finitely
generated
over $\C$.

Let us construct now the {\it homological linking form of the
deformation\/}
$$\{\ ,\ \}:T^l(M)\times T^l(M) \to \M/\OO\tag6$$
where $l$ is the middle dimension of $M$, $\ \dim M= 2l-1$ and
$\M$
denotes
the field of germs of meromorphic functions at the origin;
note that
an element of $\M$ can be represented in the form of a Laurent
series
$$f(t)=\sum_{n\ge -N}a_nt^n, \ \ a_n\in\C$$
for some non-negative integer $N$ having non-zero radius of
convergence.

Let $\M \VV$ denote $\M \otimes{\OO}\VV$ considered as a left
$\M[\pi]$-module.
Since $\VV$ is free over $\OO$ we have the following exact
sequence
$$0\to \VV\to\M \VV\to \M \VV/\VV\to 0$$
which generates the exact sequence
$$\dots \to H^{l-1}(M;\M \VV)\to H^{l-1}(M;\M
\VV/\VV)@>\delta>>H^l(M;\VV)\to
H^l(M;\M \VV)\to \dots$$
{}From this we obtain that the image of the Bockstein
homomorphism
$\delta$
is precisely the torsion  subspace $\im(\delta)= T^l$.

Finally, for $\alpha, \beta \in T^l(M)$ we can define
$$\{\alpha,\beta\}=(\delta^{-1}(\alpha)\cup\beta, [M])\tag7$$
where the cup-product is taken with respect to the natural
pairing
$$(\M \VV/\VV) \times \VV \to \M/\OO$$
determined by the Hermitian form
$\<\ ,\ \>:\VV\times \VV\to \OO$. One easily checks that
formula (7)
correctly defines a $(-1)^l$-Hermitian form
$$\{\ ,\ \}:T^l(M)\times T^l(M) \to \M/\OO$$
which is $\OO$-linear with respect to the first variable and
$\OO$-
antilinear
with respect to the second variable. Poincare duality implies
that
it is non-degenerate.

We will explain in section 2 that a linking form determines
(in a purely algebraic way) a sequence of
signature invariants
$$\sigma_1, \sigma_2,\sigma_3, \dots \in \Z.$$ These appear in
the
following
statement which is the main result of the present paper.
\proclaim {1.5. Theorem} Suppose that an analytic deformation
(3) of
a flat
Hermitian vector bundle $\E$ is given. Let $B_t$ be the
corresponding
analytic family of Atiyah--Patodi--Singer operators (2) and
let
$\eta_+,\
\eta_-$ and $\eta_0$ be the corresponding eta-invariants, cf.
(4).
Consider also the deformation of the monodromy representation
$\VV$
(cf. 1.3), corresponding to $\nabla_t$, where $t\in (-
\epsilon,\epsilon)$
and the signatures $\sigma_1, \sigma_2, \dots $
of the linking form (6).
Then the eta-invariant jumps are given by the formulae:
$$\eta_+ =\eta_0+\sum_{i\ge 1}\sigma_i, \ \ \
\ \ \eta_- =\eta_0+\sum_{i\ge 1}(-1)^i\sigma_i\tag8$$
\endproclaim

In particular, we obtain that the jumps of the eta-invariant
can be
expressed through homotopy-theoretic (even homological)
invariants.

{}From formulae (8) it follows that the {\it infinitesimal
spectral
flow}
$\eta_+-\eta_-$ is given by
$$\eta_+-\eta_-=\  2\sum_{i\ge 1}\sigma_{2i-1}\tag9$$
Note that formula (9) involves only {\it odd} signatures.
Another
nice formula which follows immediately from Theorem 1.5 is
$$\frac{1}{2}(\eta_+\ +\ \eta_-)-\eta_0=\sum_{i\ge
1}\sigma_{2i}\tag10$$
It describes the deviation of $\eta_0$ from the mean value of
$\eta_+$ and $\eta_-$ and involves only {\it even} signatures.

Similar formulae were obtained by J.Levine (cf. \cite{\Levin},
Theorem 2.3)
for jumps of the signatures of knots; the present work
actually
emerged
as the result of an attempt to understand the nature of those
jump
formulae.

The proof of Theorem 1.5 is given in section 6. The sections
2-5 are
devoted to an auxiliary material needed for the proof.

\heading 2. Linking forms and their invariants \endheading

2.1. In this section we will consider algebraic invariants of
Hermitian pairings of the form
$$\{,\}:T\otimes T\to \M/\OO\tag11$$
where $T$ is a finitely generated {\it torsion} $\OO$-module.
As
explained
in the previous section, such forms appear as linking forms
describing
deformations of monodromy representations. Very similar
algebraic
objects appear in knot theory as {\it Blanchfield pairings \/}
of
knots,
cf. \cite{\Levin}.

With any linking form (11) we will associate a {\it spectral
sequence
of quadratic forms} which will produce a set of numeriacal
invariants.

Recall that $\OO$ denotes the ring of germs of holomorphic
functions
and $\M$ denotes the field of germs of meromorphic functions.
Thus,
an
element of $\OO$ has a representation in the form of a power
series
$$f(t)=\sum_{n\ge 0}a_nt^n, \ \ a_n\in\C$$
with non-zero radius of convergency; an element of $\M$ can be
represented by a similar power series having finitely many
negative
powers. $\OO$ and $\M$ are considered together with the
involutions
induced by the complex conjugation; $t$ is assumed to be {\it
real\/},
i.e. $\overline t=t$.

The $\OO$-module $T$ in (11) viewed as a vector space over
$\C$
is finite dimensional; its $\OO$-module structure is given by
a nilpotent $\C$-linear endomorphism
$t:T\to T$ of multiplication by $t\in \OO$, with $\ t^n=0$ for
some
$n$.
The form (11) is assumed to be

(a) {\it Hermitian}, i.e $\{x,y\}=\overline{\{y,x\}}$ for all
$x,y
\in T$,
where the bar denotes the involution of $\M/\OO$ induced from
$\M$;

(b) $\OO$-{\it linear}, i.e. $\{tx,y\}=t\{x,y\}$ for $x,y \in
T$;

(c) {\it non-degenerate}, i.e. the map
$T\to \Hom_{\OO}(T,\M/\OO)$ is an isomorphism.

2.2. As an example consider the following pairing. Fix an
integer
$i\ge 1$
and a nonzero {\it real} $c$. Let $T$ be $\OO/t^i\OO$; as a
vector
space over
$\C$ it has a basis of the form $x, tx, t^2x, \dots t^{i-1}x$
where
$x$
represents the coset of $1 \in \OO$.
The pairing
$$\{,\}_{i,c}:T\otimes T\to \M/\OO\tag12$$
is given by
$$\{x,x\}_{i,c}=ct^{-i} \ \mod \OO.$$
Note that in this example only the sign of $c$ is important --
changing
$c$ to $\lambda c$ with $\lambda >0$ gives a congruent form.

2.3. Given a linking form (11), it defines the {\it scalar
form}
$$[\ ,\ ]:T\otimes T\to \C\tag13$$
where
$$[x,y]=\Res \{x,y\},$$
the residue of the meromorphic germ $\{x,y\}$. In terms of the
scalar form
$[\ ,\ ]$ one may write
$$\{x,y\}=[x,y]t^{-1}+[tx,y]t^{-2}+[t^2x,y]t^{-
3}+\dots\tag14$$
for $x,y \in T$. Thus, the scalar form contains all the
information.
It has the followings properties:

(1) $[x,y]=\overline{[y,x]}$ (i.e. it is Hermitian);

(2) the scalar form $[\ ,\ ]$ is non-degenerate;

(3) $t:T\to T$ (the multiplication by $t\in \OO$) is self-
adjoint
with
respect to the scalar form, i.e.
$$[tx,y]=[x,ty]$$
for $x,y \in T$.

2.4. Denote
$$T_i=\{x\in T; t^ix=0\}$$
for $i=0, 1, 2, \dots$. We have
$0=T_0\subset T_1 \subset T_2 \dots $
and $T=T_N=T_{\infty}$
for large $N$.

Note that
$T_i\supset tT_{i+1},\ $ and the natural inclusion map
$T_i\to T_{i+1}$ induces an {\it inclusion} $T_i/tT_{i+1}\to
T_{i+1}/tT_{i+2}.$
Denote $$V_i=T_i/tT_{i+1};$$ for large $i$ the space $V_i$ is
equal
to
$V_{\infty}=T/tT.$ Thus we have the sequence of vector spaces
$$0=V_0\subset V_1 \subset V_2 \ \dots \subset V_{\infty}$$
On every $V_i$ there is defined a Hermitian form
$$l_i:V_i \times V_i \to \C$$
where
$$l_i(x,y)=[t^{i-1}x,y]$$
for $x,y \in T_i;$ one easily checks that this formula
correctly
defines a form on $V_i=T_i/tT_{i+1}$.

\proclaim{2.5. Lemma} The annihilator of the form $l_i:V_i
\times
V_i \to \C$
is equal to $\ V_{i-1}$, i.e.
$$\ker (l_i) = V_{i-1}.$$
\endproclaim
\demo{Proof} Using non-degeneracy of the scalar form $[\ ,\ ]$
one
first
checks that
$T_i^{\bot}=t^iT.$
If $x\in T_i $ and $[t^{i-1}x,y]=0$ for any $y\in T_i$ then
$t^{i-1}x\in T_i^{\bot}=t^iT$ and thus
$t^{i-1}x=t^iz$ where $z\in T_{i+1}$. We obtain that
$x=tz +u$
with $u\in T_{i-1}$. These arguments show that
$V_{i-1}\supset \ker (l_i)$. The other inclusion is obvious.
$\square$
\enddemo

2.6. Thus we obtain that any linking form
$$\{\ ,\ \}:T\otimes T\to \M/\OO$$
defines an algebraic object which
we will call {\it
a spectral sequence of quadratic forms}
(because of its similarity to spectral sequences).
It consists of a flag of finite dimensional vector spaces
$$0=V_0\subset V_1 \subset V_2 \ \dots \subset
V_{\infty}=T/tT$$
(which actually stablize, $V_i=V_{\infty}$ for large $i$)
supplied with a sequence of Hermitian forms
$$l_i:V_i \times V_i \to \C, \ \ i=0,1,2,\dots$$
such that

(1) $l_i$ vanishes identically for large $i$;

(2) $V_{i-1}= \ker (l_i)$.

2.7. Using the spectral sequence of quadratic forms
associated to a linking form one may construct a set of
numerical
invariants of linking forms.
For any integer $i\ge 1$ let $n_i^{+}$ and $n_i^{-}$ denote
the number of positive and negative squares appearing in the
diagonalization of the Hermitian form
$l_i:V_i \times V_i \to \C$. Let $\sigma_i$ denote the
difference
$$\sigma_i= n_i^{+}-n_i^{-};$$
it is the signature of $l_i$. The numbers
$$n_i^+, n_i^-,\sigma_i$$
are obviously {\it invariants of the linking form}. The
invariants
$\sigma_i$ will be the most important for the sequel; we will
call
them
{\it signatures \/} of the form (11).

2.8. There is yet another spectral sequence of quadratic forms
associated
with any linking form (11), which is in fact more useful.

Denote $W_i=t^{i-1}T_i$ for $i=1,2,\dots$. Thus, $W_1=T_1$ and
we
have
a decreasing finite filtration
$$W_1\supset W_2\supset W_3\supset\dots \supset W_\infty =0.$$
For any integer $i\ge 1$ define a Hermitian form
$$\lambda_i: W_i\times W_i\to \C$$
by $\lambda_i(x,y)=\Res \{a,y\}$ where $a\in T_i$ is an
element
with $t^{i-1}a=x$.
Clearly, the form $\lambda_i$ is correctly defined and is
Hermitian.

The following statement is similar to Lemma 2.5; it states
that the
set
of forms $\lambda_i$ form a {\it spectral sequence of
Hermitian
forms}.
Note that, the spectral sequence formed by $W_i$'s, grows in
the
opposite
direction compared with the spectral sequence of 2.6.

\proclaim{2.9. Lemma} The annihilator of the form
$\lambda_i: W_i\times W_i\to \C$ is equal to $W_{i+1}$ and the
induced
nondegenerate form  on $W_i/W_{i+1}$ has exactly $n_i^+$
positive
squares
and $n_i^-$ negative squares, where the numbers $n_i^+$ and
$n_i^-$
are
defined in 2.7. Thus, the signature $\sigma_i$ can be also
computed
as
the signature of the form $\lambda_i$.
\endproclaim
\demo{Proof} By Lemma 2.5, $l_i$ induces a nondegenerate form
on
$V_i/V_{i-1}$; let us denote this induced form $\overline
l_i$. On
the other
hand, observe that
$\lambda_i(x,y)=\Res \{a,y\}=0$ if $x\in W_i$ and $y\in
W_{i+1}$.
Thus, $\lambda_i$ induces a form $\overline\lambda_i$
on the factor $W_i/W_{i+1}$. We claim that there is an
isomorphism
$$\alpha_i: W_i/W_{i+1}\ \to \ V_i/V_{i-1}$$
which intertwins between $\overline l_i$ and
$\overline\lambda_i$.
This
would obviously imply the statement of the Lemma.

If $x\in W_i$, represent $x$ as $t^{i-1}a$ for some $a\in T_i$
and
define $\alpha_i(x)$ to be equal to the coset of $a$ in
$V_i/V_{i-1}=T_i/(T_{i-1}+tT_{i+1})$. One easily checks that
this map is correctly defined and has the properties mentioned
above.
$\square$
\enddemo

2.10. As an example consider the linking form (12)
$$\{,\}_{i,c}:T\otimes T\to \M/\OO,$$
where $T=\OO/t^i\OO$.
In this case we obtain that
$V_j= 0$ for $j<i$ and $V_j=\C$ for $j\ge i$.
All forms $l_j$ with $j\ne i$ vanish; the $i$-th form $l_i$
has
signature equal to the sign of the number $c$.

Since the invariants $n_i^+, n_i^-,\sigma_i$ are {\it
additive} we
obtain
that:

\proclaim{2.11. Corollary} Given a linking form (11)
which is represented as orthogonal
sum of finitely many forms of type $\{,\}_{i,c}$ (with
different $i$
and $c$),
the number $n_i^+$ is equal to the number of summands
of type
$\{,\}_{i,c}$ with $c$ positive and the number $n_i^-$
is equal to the number of summands $\{,\}_{i,c}$ with $c$
negative
in the above decomposition.
\endproclaim

2.12. It is easy to show that {\it any linking form (11) is
diagonalizable,
i.e. it is congruent to a direct
sum of forms of the type $\{,\}_{i,c}$}; we will not use this
fact
in the present paper and thus will leave it without proof.
The uniqueness of this
orthogonal decomposition follows from the previous arguments.
Thus the numbers $n_i^+$ and $n_i^-$ {\it determine the type
of the
form
(11)}.

A linking form
$$\{\ ,\ \}:T\otimes T\to \M/\OO$$
will be called {\it hyperbolic\/} if the $\OO$-module $T$ can
be
represented as a direct sum $T=A\oplus B$ such that the
restrictions
of the form $\{\ ,\ \}$ on $A$ and on $B$ vanish: $\{x,y\}=0$
if
either
$x,y\in A$ or $x,y\in B$ (this can be expressed by saying that
$A$
and
$B$ are Lagrangian direct summands).

\proclaim{2.13. Lemma} All signatures $\sigma_i,\ i\ge 1$
of a hyperbolic linking form vanish.
\endproclaim
\demo{Proof} Suppose $T=A\oplus B$ where $A$ and $B$ are
Lagrangian
direct summands (over $\OO$). Then for any integer $i\ge 1$
the vector space $T_i$ (defined as in 2.4)
is also a direct sum $T_i=A_i\oplus B_i$ of vector spaces
defined in the similar way by
$A$ and $B$ respectively; thus the
vector space $V_i=V_i(T)$ is also given as a direct sum
$V_i(T)=V_i(A)\oplus V_i(B)$. By Lemma 2.5 the pairing $l_i$
induces a {\it non-degenerate\/} pairing $\tilde l_i$ on
$V_i(T)/V_{i-1}(T)$
and the signature of $\tilde l_i$ is equal to $\sigma_i$. Thus
we
obtain
that $V_i(T)/V_{i-1}(T)$ is a direct sum
$$V_i(A)/V_{i-1}(A)\oplus V_i(B)/V_{i-1}(B)$$
and the the form $\tilde l_i$ vanishes on
$V_i(A)/V_{i-1}(A)$ and on $V_i(B)/V_{i-1}(B)$. This implies
that
$\sigma_i=0$ for all $i\ge 1$. $\square$
\enddemo

2.14. We also mention the closely related notion of {\it
metabolic}
form.
By the definition, a form
$\{\ ,\ \}$ is metabolic if there is a submodule $A\subset T$,
of
half the
dimension (as vector space over $\C$) of $T$, such that $\{\
,\ \}$,
restricted to $A$, vanishes. (We say $A$ is a {\it
Lagrangian}).
When
$\{\ ,\ \}$ is metabolic the individual signatures are not
necessarily zero
but we do have the property $\sum_{i\ge 1}\sigma_{2i-1}=0$. In
fact,
this equation is a necessary and sufficient condition for $\{\
,\
\}$
to be metabolic. We will not prove these statements here since
they
will
not be used in the present work. It is interesting to compare
this
to
formula (9) above. The relations between hyperbolic and
metabolic
forms and
related signature invariants are discussed more fully, in a
special
case,
in \cite{\Levin}.

2.15. As a concluding remark let us note that the
study of skew-Hermitian forms (11) (i.e. forms satisfying
$\{x,y\}=-\overline{\{y,x\}}$ for all $x,y\in T$ together with
(b)
and (c)
of subsection 2.1) can be automatically reduced
to the case of Hermitian forms discussed above by multiplying
the form by $i=\sqrt{-1}$. Thus we {\it define\/}
signatures $\sigma_j,\ \ j\ge 1$ of a skew-Hermitian
linking form (11) as the corresponding signatures of the form
$i\{\ ,\ \}$.

\heading 3. Jumps of the eta-invariant and signatures of the
linking form determined by deformation of a self-adjoint
operator\endheading

We are going to establish that
an analytic deformation
of a self-adjoint elliptic operator defines a linking form
of the type studied in the previous section. We will call this
linking
form analytic in order to distinguish it from the algebraic
linking
form
defined in 1.4. We will prove
(it will be the main result of this section) that the jumps of
the eta-invariant of the family of operators
can be expressed through a combination of signatures
associated to the analytic linking form. We will also compute
explicitly
the corresponding spectral sequence of quadratic forms in
terms of
Taylor
expansion of the family.

The above mentioned linking form
is constructed by studying the action of the family of the
operators
on germs of analytic curves of sections of a vector bundle.
The idea of considering the family of operators as a single
operator
acting on the space
of curves is actually
the principal technical novelty of the present paper.

We start this section by defining
precisely the analytic curves we are going to use.

\subheading{3.1} First we recall some standard definitions.

Let $\Omega$ be an open subset of $\C$ and let $V$ be
a complex topological vector space. A function $f:\Omega\to V$
is
said to be
{\it weakly holomorphic in $\Omega$} if $vf$ is holomorphic in
the
ordinary sense
for every continuous linear functional $v$ on $V$. The
function
$f:\Omega\to V$ is called {\it strongly holomorphic in
$\Omega$} if
the limit
$$\lim_{w\to z}\frac{f(w)-f(z)}{w-z}$$
exists (in the topology of $V$) for every $z\in\Omega$. It is
known
that
the above two notions of analyticity actually
coincide if $V$ is a Frechet space, cf. \cite{\Rud}, Chapter
3.

A function $f:(a,b)\to V$ defined on a real interval $(a,b)$
with
values in a
Frechet space $V$ is called
{\it analytic (or real analytic or holomorphic)} if it is a
restriction
of an analytic function $\Omega\to V$ defined in a
neighbourhood
$\Omega\subset\C$ of the interval $(a,b)$.

\subheading{3.2}
We will mainly consider analytic curves in spaces of smooth
sections of vector bundles.
Let $M$ be a compact $C^\infty$ Riemannian
manifold (possibly with
boundary) and let $\E$ be a Hermitian vector bundle over $M$.
For
any
integer $k$ symbol $\H_k(\E)$ will denote the corresponding
Sobolev
space
(defined as in Chapter 9 of \cite{\Pal}). Recall that the
Sobolev
spaces
$\H_k(\E)$ with $k\in\Z$ form a chain of Hilbertian spaces (in
the
terminology
of \cite{\Pal}), which, in particular, means that $\H_k(\E)$
is
embedded
into $\H_l(\E)$ for $k>l$ (as a topological vector space)
and the intersection of all the spaces $\H_k(\E)$
coincides with $\H_{\infty}(\E)=C^\infty(M)$.

\subheading{3.3. Definition}
Let $f:(a,b)\to C^\infty(M)$ be a curve of smooth sections;
we will say that $f$ is {\it
analytic} if for any integer $k$ the curve $f$ represents a
(real) analytic curve
considered as a curve in the Sobolev space $\H_k(\E)$.

Note that any curve $f:(a,b)\to C^\infty(M)$ which is analytic
in
the
sense of section 3.1 (i.e. by viewing $C^\infty(M)$ as a
Frechet
space)
will be obviously analytic in the sense of Definition 3.3. The
converse
is also true, although the proof of this fact
is not elementary; we are grateful to
V.Matsaev for explaining this to us. The proof suggested by
V.Matsaev
uses interpolation theory of Hilbert spaces. Since we wish to
avoid these analytic subtleties, and since the definition 3.3
is
the most convenient and entirely sufficient for our purposes,
we will accept it and will never use
the equivalence of the above two definitions in the present
paper.

Suppose that $\E$ and $\F$ are two Hermitian vector bundles
over
$M$. Then
any differential operator $D:C^\infty(\E)\to C^\infty(\F)$ of
order
$\ell$
defines a bounded linear map of Sobolev spaces $\H_k(\E)\to
\H_{k-
\ell}(\F)$
(where $k\ge \ell$)
and thus $D$ maps analytic
curves in $C^\infty(\E)$ into analytic curves in
$C^\infty(\F)$.

\subheading{3.4} We will give now
definition of analyticity for families
of linear differential operators.

Suppose that $\E$ and $\F$ are two vector bundles over the
manifold
$M$ and
$D_t\in \Diff_\ell(\E,\F)$ is a family of linear differential
operators of
order $\ell$, depending on a real parameter $t\in (a,b)$.
Let $J^\ell(\E)$ denote the jet bundle of order
$\ell$, cf. \cite{\Pal}, chapter IV, \S 2. Then by Theorem 1
on page
61
of \cite {\Pal}, the set $\Diff_\ell(\E,\F)$ can be identified
with
$C^\infty(\Hom(J^\ell(\E),\F))$. The latter is the set of
smooth
sections
of a vector bundle; therefore we can consider analytic curves
in
this space of
sections using the definition of analyticity given in 3.3.

We accept the following definition:
a curve of linear differential operators $(a,b)\to
\Diff_\ell(\E,\F)$
is called {\it (real) analytic}
iff the corresponding curve of sections of the bundle
$\Hom(J^\ell(\E),\F)$ is analytic.

The main property of analytic families of operators $D_t$,
which we
will
constantly use, consists of the following: for any integer
$k\ge\ell$
the family of bounded
linear operators $D_t: \H_k(\E)\to \H_{k-\ell}(\F)$ depends
analytically
on the parameter $t$ (i.e. defines an analytic curve in the
Banach
space
of bounded linear operators $\H_k(\E)\to \H_{k-\ell}(\F)$
with the operator norm).

{}From the above remark it follows that {\it
if $f:(a,b)\to C^\infty(\E)$ is an
analytic curve of smooth sections and if $D:(a,b)\to
\Diff_\ell(\E,\F)$
is an analytic curve of linear differential operators then the
\lq\lq evaluation curve\rq\rq $t\mapsto D_t(f_t)$ is also
analytic.}

We will formulate now a few simple lemmas which will be used
later.
Roughly speaking, they represent different converses of the
statement of the
previous paragraph.

\proclaim{3.5. Lemma} Let $\E$ be a Hermitian vector bundle
over a
compact Riemannian manifold $M$ without boundary and let
$D_t\in \Diff_\ell(\E,\E)$
be an analytic (in the sense of 3.4) family of elliptic self-
adjoint
operators of order
$\ell>0$ defined for $t\in (a,b)$. Suppose that $\ker D_t=0$
for all
$t\in (a,b)$. If $\phi, \psi :(a,b)\to C^\infty(\E)$
are two curves such that $D_t(\phi(t))=\psi(t)$ for any $t\in
(a,b)$
and
the curve $\psi$ is analytic (in the sense of definition 3.3),
then
the
curve $\phi$ is also analytic.
\endproclaim
\demo{Proof} Fix an integer $k$.
Since $\ker D_t=0$, the operator $D_t$
defines a linear homeomorphism $D_t: \H_{k+\ell}(\E)\to
\H_k(\E)$
(by the
open mapping theorem, cf. \cite{\Rud}, p.47)
which depends analytically on $t$. Thus it follows that
$\phi(t)=D_t^{-1}(\psi(t))$
is an analytic curve in the Sobolev space $\H_{k+l}(\E)$.
Since this is true for any $k$, the statement follows.
$\square$
\enddemo
\proclaim{3.6. Lemma} Let $\E$ be a Hermitian vector bundle
over a
compact Riemannian manifold $M$ without boundary and let
$D_t\in \Diff_\ell(\E,\E)$
be an analytic (in the sense of 3.4) family of elliptic self-
adjoint
operators of order
$\ell>0$ defined for $t\in (a,b)$.
Suppose that $\phi, \psi :(a,b)\to C^\infty(\E)$
are two curves such that $D_t(\phi(t))=\psi(t)$ for any $t\in
(a,b)$
and
it is known that the curve $\psi$ is analytic in the sense of
definition 3.3, while
the curve $\phi$ is analytic in a weaker sense -
as a curve in the Hilbert space $\H_0(\E)=L^2(\E)$. Then the
curve $\phi$ is analytic in the sense of definition 3.3 as
well.
\endproclaim
\demo{Proof} Choose a point $t_0\in (a,b)$ and an integer
$k\ge 0$.
It is enough to prove analyticity of the curve
$\phi :(t_0-\delta,t_0+\delta)\to \H_k(\E)$
for some small $\delta>0$ (the restriction of the original
curve
$\phi$
onto a neighbourhood of $t_0$, considered as a curve in the
Sobolev
space
$\H_k(\E)$).

Let $\pi$ denote the orthogonal projection of $\H_0(\E)$ onto
$\ker(D_{t_0})\subset \H_\infty(\E)$. The operator
$$D_t+\pi:\quad \H_k(\E)\to\H_{k-\ell}(\E)$$
is continuous, analytically depends on the parameter $t$, and
is
invertible
for $t=t_0$. Thus it is invertible for $t\in (t_0-
\delta,t_0+\delta)$
for some $\delta>0$. We have
$$(D_t+\pi)(\phi(t))=\psi(t)+\pi(\phi(t))$$

We claim that the right hand side of this equation is a curve
analytic
in the Sobolev space
$\H_{k-\ell}(\E)$. In fact, the first summand $\psi(t)$ is
analytic
in any
Sobolev space by the
assumption, while the second summand $\pi(\phi(t))$ belongs to
a
finite
dimensional subspace $\ker D_{t_0}$, and it is given that
it is analytic as a curve in Hilbert space
$L^2(\E)=\H_0(\E)$. Since all linear topologies on a finite
dimensional
vector space are equivalent, we conclude that the curve
$\pi(\phi(t))$
is analytic as a curve in $\H_{k-\ell}(\E)$.

Combining the remarks of the two previous paragraphs, we
obtain that
the
curve $\phi :(t_0-\delta,t_0+\delta)\to \H_k(\E)$ is analytic.
$\square$.
\enddemo
\subheading{3.7} Suppose again that $\E$ is a Hermitian vector
bundle over a
compact Riemannian manifold $M$ without boundary and
$D_t\in \Diff_\ell(\E,\E)$ is an analytic (in the sense of
3.4)
family of elliptic self-adjoint operators of order $\ell>0$
defined for $t\in I=(a,b)$. In this situation there exists a
parametrized
spectral decomposition (cf. \cite{\Kato}, Theorem 3.9, Chapter
VII,
\S 3)
which consists in a sequence of analytic (in
$L^2(\E)=\H_0(\E)$)
curves $\phi_n(t)$ and a sequence of analytic real valued
functions
$\mu_n(t)$ (defined for all $t\in I$) such that $\mu_n(t)$
represent
all
the repeated eigenvalues of $D_t$ and $\phi_n(t)$ form a
complete
orthonormal
family of the associated eigenvectors of $D_t$ acting on
Hilbert
space
$L^2(\E)=\H_0(\E)$. By the regularity theorem for elliptic
operators, the
curves $\phi_n(t)$ actually belong to $C^\infty(\E)$.

We claim now that {\it the curves of eigenfunctions
$\phi_n(t)$, which appear in the parametrized spectral
decomposition, are analytic in the sense of
Definition 3.3, i.e. as curves in any Sobolev space
$\H_k(\E)$}.
In fact, it is
enough to apply Lemma 3.6 to the equation
$$(D_t-\mu_n(t))(\phi_n(t))=0$$
and observe that the operator $D_t-\mu_n(t)$ depends
analytically on
$t$
while the curve on the right (the zero curve) is analytic.

\subheading{3.8. Construction of the linking form}
Let $\E$ be a Hermitian vector bundle over a closed Riemannian
manifold $M$.
Symbol $\OO C^\infty(\E)$ will denote
the set of germs at $t=0$ of all analytic curves $(-
\epsilon,\epsilon)\to
C^\infty(\E)$ in the sense of
definition 3.3. It is a left module over the ring $\OO$ of
germs of
analytic curves on the complex plane (via the pointwise
multiplication),
which has already appeared in \S 2. There is also a
\lq\lq scalar product\rq\rq
$$(\ ,\ ):\OO C^\infty(\E)\times \OO C^\infty(\E)\to
\OO\tag15$$
(the pointwise scalar product of curves of sections),
which is $\OO$-linear with respect to the first variable and
skew-
linear
with respect to the second variable.

Suppose that $D_t\in\Diff_\ell(\E,\E)$ is an analytic curve of
elliptic
self-adjoint differential operators of order $\ell >0$ defined
for
$t\in (-\epsilon,\epsilon)$, cf. 3.4. Then it defines the
following
single map
$$\tilde D:\quad \OO C^\infty(\E)\to \OO C^\infty(\E)\tag16$$
where for $\alpha\in \OO C^\infty(\E)$ the germ $\tilde
D(\alpha)$
represents the curve $t\mapsto D_t(\alpha(t))$. It is clear
that
$\tilde D$
is an $\OO$-homomorphism. Consider the image of $\tilde D$ and
a
larger
$\OO$-submodule ${\frak Cl}(\im(\tilde D))\subset\OO
C^\infty(\E)$
consisting of germs $\alpha$ with
the property that $t^k\alpha$ belongs to $\im(\tilde D)$ for
some
$k>0$;
here $t^k$ denotes the element of the ring $\OO$ represented
by the
curve
$t\mapsto t^k$. Now define
$$T={\frak Cl}(\im(\tilde D))/\im(\tilde D);\tag17$$
it is a module over $\OO$.

Let $\M$ denote the field of fractions of $\OO$; in other
words,
$\M$
is the field of germs at $0$ of
meromorphic curves on $\C$. Let us define the {\it
linking pairing}
$$\{\ ,\}:T\times T\to \M/\OO.\tag18$$
If $\alpha,\beta$ are two given elements of $T$, represent
$\alpha$
and
$\beta$ by germs of curves $f$ and $g$ in
$\OO C^\infty(\E)$ correspondingly;
then $t^kf=\tilde D(h)$ for some $k>0$ and
$h\in\OO C^\infty(\E)$. Now we define
$$\{\alpha,\beta\}=t^{-k}(h,g)\quad\in\quad\M/\OO\tag19$$
One easily checks that the definition is correct. We will also
refer
to (18) as the {\it analytic linking form\/}
associated with the deformation $D_t$.

\proclaim{3.9. Theorem} (1) The analytic linking pairing (18)
(constructed out of
an analytic family of elliptic self-adjoint
operators $D_t:C^\infty(\E)\to C^\infty(\E)$, where $-
\epsilon<t<\epsilon$,
acting on sections of a Hermitian vector bundle $\E$ over a
closed
Riemannian manifold $M$),
is Hermitian, $\OO$-linear with
respect to the first variable and non-degenerate
and the $\OO$-module $T$ is finite-dimensional over $\C$;
(2) Let $\eta(D_t)$ denote the eta-invariant of the operator
$D_t$
(cf. \cite{\Gilkey}) and let $\eta_{\pm}$ denote the limits
$$\eta_{\pm}=\lim_{t\to {\pm}0}\eta(D_t).\tag20$$
Then the following jump formulae hold:
$$\eta_+ = \eta_0 + \sum_{i\ge1} \sigma_i,\ \ \
\eta_- = \eta_0 +\sum _{i\ge 1} (-1)^i \sigma_i;\tag21$$
here $\eta_0$ denotes $\eta(D_0)$ and $\{\sigma_i\}$ denote
the
signatures
of the linking pairing $\{\ ,\ \}$, introduced in \S 2.
\endproclaim
\demo{Proof} Consider the parametrized spectral decomposition
$\phi_n(t),\mu_n(t)$ where $n=1,2,\dots$,
and $t\in I=(-\epsilon,\epsilon)$,
of the analytic self-adjoint elliptic family $D_t$,
mentioned in subsection 3.7. For every value of $t\in I$ the
vectors
$\phi_n(t)$ form a complete orthonormal system of eigenvectors
of
$D_t$
with eigenvalues $\mu_n(t)$ and we know that $\mu_n(t)$ are
analytic
functions of $t\in I$ and the curves of sections $\phi_n(t)$
are
analytic
in the sense of definition 3.3, as shown in 3.7.

{}From ellipticity of $D_t$ it follows that there exist only
finitely
many numbers $n$ such that $\mu_n(0)=0$. We can suppose that
the
numeration
of the eigenfunctions and the corresponding eigenvalues has
been
arranged
so that
\roster
\item  $\mu_n(0)=0$ while $\mu_n(t)$ are not identically zero
for $n=1,2,\dots, N$;
\item $\mu_n(t)\equiv 0$ for $n=N+1, N+2,\dots, N_1$;
\item $\mu_n(0)\ne 0$ for all $n>N_1$.
\endroster

The proof of the Theorem will be based on the following
statement:

\proclaim{Claim} For any curve $\psi:(-\epsilon,\epsilon)\to
C^\infty(\E)$,
which is analytic in the sense of Definition 3.3 and satisfies
$$(\phi_n(t),\psi(t))=0\quad\text{for}\quad n=1,2,\dots, N_1$$
(where we use $L^2(\E)=\H_0(\E)$-scalar product), there exists
a curve $\phi:(-\epsilon,\epsilon)\to C^\infty(\E)$ satisfying
\roster
\item $D_t(\phi(t))=\psi(t)\quad\text{for all}\quad t\in (-
\epsilon,\epsilon)$;
\item $(\phi_n(t),\phi(t))=0\quad\text{for}\quad n=1, 2,
\dots,
N_1$;
\endroster
and the curve $\phi:(-\epsilon,\epsilon)\to C^\infty(\E)$ is
analytic
in the sense of 3.3.
\endproclaim
Such curve $\phi$ is actually unique but we will not need this
fact.
\demo{Proof of the Claim} For any $t\in (-\epsilon,\epsilon)$
let $\pi_t: \H_0(\E)\to\H_\infty(\E)$ denote
the orthogonal projection onto the finite
dimensional subspace generated by
$\phi_1(t),\dots,\phi_{N_1}(t)$.
Represent $\pi_t=\pi_t^\prime+\pi_t^{\prime\prime}$ where
$\pi_t^\prime$
is the projection onto the subspace generated by the curves
$\phi_n(t)$, with $1\le n\le N$, and
$\pi_t^{\prime\prime}$ is the projection onto
the subspace generated by the curves $\phi_n(t)$ with
$N+1\le n\le N_1$. These operators define naturally operators
acting
on all Sobolev spaces.
{}From the fact that the curves $\phi_n(t)$ are analytic in the
sense
of 3.3
(cf. 3.7) it follows that the operators
$$\pi_t^\prime,\ \pi_t^{\prime\prime}: \H_k(\E)\to
\H_\infty(\E)$$
analytically depend on $t$ for any integer $k$.

Consider the operators
$$(D_t+\pi_t):\H_k(\E)\to \H_{k-\ell}(\E);$$
they form an analytic family of linear homeomorphisms
(i.e. an analytic curve in
the space of bounded linear operators $\H_k(\E)\to \H_{k-
\ell}(\E)$
with
the operator norm topology). It follows that the equation
$$(D_t+\pi_t)(\xi(t))=\psi(t)$$
admits a unique solution $\xi(t)$ which lies in $C^\infty(\E)$
and
is
analytic in the sense of definition 3.3. Since
$\pi_t(\psi(t))=0$
and
$\pi_t^{\prime\prime}(D_t(\xi(t)))=0$ we obtain
$$\pi_t^\prime(D_t(\xi(t)))+
\pi_t^\prime(\xi(t))=0\quad\text{and}\quad
\pi_t^{\prime\prime}(\xi(t))=0.$$
Moreover, since  $\pi_t^\prime(\xi(t))$ is an eigenvector of
$D_t$
with
eigenvalue -1 and on the other hand it is a sum eigenvectors
each of
whose eigenvalues has constant term 0, we obtain
$\pi_t^\prime(\xi(t))=0$.
Thus, the curve $\phi(t)=\xi(t)$
$D_t\phi=\psi,\quad \pi_t\phi=0$. This completes the proof
of the Claim.
\enddemo

Now we may easily finish the proof of Theorem 3.9.
Given a germ of an analytic curve $f\in \OO C^\infty(\E)$ we
can
find
(using the Claim) a germ $g\in \OO C^\infty(\E)$ with
$$f-\pi_tf=\tilde D(g),\quad \pi_t(g)=0.$$
Here we use the notation introduced in the proof of the Claim.
Since the the projector $\pi_t^{\prime\prime}$ vanishes on
${\frak Cl}(\im(\tilde D))$ we obtain that:

{\it
(i) a germ $f\in\OO C^\infty(\E)$ belongs to $\im(\tilde D)$
if and
only if
$\pi_t^\prime(f)\in\im(\tilde D)$ and
$\pi_t^{\prime\prime}(f)=0$;

(ii) a germ $f\in\OO C^\infty(\E)$ belongs to ${\frak
Cl}(\im(\tilde
D))$
if and only if
$\pi_t^\prime(f)\in{\frak Cl}(\im(\tilde D))$ and
$\pi_t^{\prime\prime}(f)=0$.}

Now we see that the germs of eigenfunctions
$\phi_1(t),\phi_2(t),\dots
\phi_N(t)$ generate $T$ as a module over $\OO$. The full set
of
relations
of $T$ is given by $\mu_n(t)\phi_n(t)=0$, where $n=1,2,\dots,
N$.

Computing the linking form on the elements $\phi_n\in T$,
where
$1\le n\le N$, representing
the eigenfunctions, we obtain
$$
\{\phi_i,\phi_j\}=\cases 0, &\text{for $i\ne j$}\\
(\mu_i(t))^{-1} \ \mod \OO,&\text{for $i=j$}.\endcases
\tag22
$$

We also obtain that {\it the invariant $n^+_i$ of the linking
form
(18) (which was defined in \S 2) is equal to the number of
eigenvalues
$\mu_n(t)$ having the form $\mu_n(t)=t^i\overline \mu_n(t)$
where
$\overline \mu_n(0)$ is positive. Similarly,
the invariant $n^-_i$ of the linking form
is equal to the number of eigenvalues
$\mu_n(t)$ having the form $\mu_n(t)=t^i\overline \mu_n(t)$
where
$\overline \mu_n(0)$ is negative.}

Since $\sigma_i=n^+_i-n^-_i$,
we obtain that the jump of the eta-invariant
$\eta_+ - \eta_0$ is given by $ \sum_{i\ge1} \sigma_i$.
The second formula follows similarly. $\square$
\enddemo

3.10. We are going now to compute explicitly the spectral
sequence
of
quadratic forms (described in 2.8) of the linking form of
self-
adjoint
analytic family $D_t$, constructed in 3.8. We will show that
this
spectral
sequence of quadratic forms can be expressed in terms of the
kernel of the undeformed operator $\ker (D_0)$ (the space of
\lq\lq harmonic forms\rq\rq)
and the pairings on it
given by the terms of the Taylor expansion of the family
$$D_t=D_0 + tD_1 + t^2D_2 +\dots$$
We want to emphasize that each form of the spectral sequence
uses
only
finitely many derivatives $D_i$.

First, we are going to identify the initial term
$W_1=T_1=\{\alpha\in T;
t\alpha=0\}$ of this spectral sequence, cf. 2.8.
According to the definition of subsection 3.8,
any element $\alpha\in T_1$ can be represented by an analytic
germ $\alpha_t\in C^\infty(\E)$ such that
$$t\alpha_t = D_t(\beta_t)$$
for some analytic germ $\beta_t\in C^\infty(\E)$. The curve
$\beta_t$
determines $\alpha$; adding to $\beta_t$ a curve of the form
$t\gamma_t$ does not change the class of $\alpha$ in $T$.
Thus, we obtain that
the initial term $\beta_0$ of the Taylor expansion
$$\beta_t=\beta_0+t\beta_1+t^2\beta_2+\dots$$
determines the class $\alpha\in W_1$. The section $\beta_0$
must
satisfy
$D_0(\beta_0)=0$. Let $\Sigma\subset \ker(D_0)$ denote the set
of
all
$s\in C^\infty(\E)$ satisfying $D_i(s)=0$ for every $i\ge 0$.
Thus
we
obtain that:

{\it the initial term of the spectral sequence of quadratic
forms
associated with the linking form of an analytic self-adjoint
family,
is given by}
$$W_1\ =\ \ker (D_0)/\Sigma.$$

Let us now compute the first form $\lambda_1: W_1\times W_1\to
\C$,
defined in
2.8. Suppose that $\beta_0,\beta^\prime_0\in\ker (D_0)/X$ are
two
elements
representing $\alpha,\alpha^\prime\in T_1=W_1$
correspondingly.
Then, we have
$t\alpha_t=D_t(\beta_0)$ and
$t\alpha_t^\prime=D_t(\beta_0^\prime)$
and,
combining our definitions, we get
$$\lambda_1(\beta_0,\beta_0^\prime)=\Res\{\alpha,\alpha^\prime
\}=
\Res t^{-1}(\beta_t,\alpha_t^\prime)=
(\beta_0,D_1(\beta_0^\prime))$$
Thus, we obtain the folowing statement, which in the case of a
particular
operators on three dimensional manifolds was established by
P.Kirk and E.Klassen \cite{\Kirk}:

\proclaim{3.11. Corollary} The first Hermitian form
$\lambda_1: W_1\times W_1\to \C$ can be identified with the
form on
$\ker (D_0)/\Sigma$ induced by the first derivative $D_1$. In
particular, the
first signature $\sigma_1$, as well as the invariants $n_1^+$
and
$n_1^-$,
can be computed as the corresponding invariants of the
Hermitian
form
on $\ker(D_0)$ given by $(x,y)\mapsto (D_1x,y)$.
\endproclaim

By the construction of 2.8, the annihilator of $\lambda_1$ is
$W_2$,
and
there is a form $\lambda_2$ defined on $W_2$ with annihilator
$W_3$,
and so
on. We will desribe all these forms $\lambda_i$ as follows.

3.12. By 2.8, $W_i=t^{i-1}T_i\subset W_1$. An element $\alpha$
of
$T_i$
can be represented by a curve $\alpha_t\in C^\infty(\E)$ such
that
$t^i\alpha_t=D_t(\beta_t)$ for some analytic curve $\beta_t\in
C^\infty(\E)$.
We observe that the first $i$ coefficients of the Taylor
expansion
$\beta_t=\beta_0+t\beta_1+t^2\beta_2+\dots$ determine the
class of
$\alpha\in T_i$. In other words, an element of $T_i$ can be
described by
a polynomial curve
$$\beta_t=\beta_0+t\beta_1+\dots +t^{i-1}\beta_{i-1}$$
of degree $i-1$ satisfying the following system of equations
$$
\left\{
\aligned D_0\beta_0&=0\\
D_1\beta_0+D_0\beta_1&=0\\
\dots\quad\dots&\\
D_{i-1}\beta_0+D_{i-2}\beta_1+\dots +D_0\beta_{i-1}&=0
\endaligned\right.\tag23_i
$$

Thus, we obtain that {\it an element of $W_i=t^{i-1}T_i$ can
be
identified with the
set of all $\beta_0\in\ker(D_0)/\Sigma$, such that the system
of
equations
$(23_i)$ admits a solution with given $\beta_0$.}

Suppose now that $\beta_0,\beta_0^\prime\in W_i$ are two such
elements.
We want to compute their product
$\lambda_i(\beta_0,\beta_0^\prime)$.
Denote $\beta_t= \beta_0+t\beta_1+\dots +t^{i-1}\beta_{i-1}$
where $\beta_1,\beta_2,\dots\beta_{i-1}$ form a solution of
$(23_i)$
with
given $\beta_0$.
Then $D_t\beta_t=t^i\alpha_t$ for some analytic curve
$\alpha_t$.
Define similarly $\beta_t^\prime$ and $\alpha_t^\prime$
According to our definitions given in 2.8, 3.8, we obtain
$$
\aligned
\lambda_i(\beta_0,\beta_0^\prime)=&\Res \{\alpha_t,t^{i-
1}\alpha_t^\prime\}\\
=&\Res \{t^{i-1}\alpha_t,\alpha_t^\prime\}\\
=&\Res t^{-1}(\beta_t,\alpha_t^\prime)\\
=&(\beta_0,\alpha_0^\prime)\\
=&(\beta_0,D_i\beta_0^\prime+D_{i-1}\beta_1^\prime+\dots+
D_1\beta_{i-1}^\prime)\endaligned
$$
Thus we have established
\proclaim{3.13. Corollary} The $i$-th term $W_i\subset
\ker(D_0)/\Sigma$
of the spectral sequence of quadratic forms determined by the
linking form
(18) consists of all $\beta_0$ such that the system $(23_i)$
has a
solution
with given $\beta_0$. If $\beta_0,\beta_0^\prime\in W_i$ are
two
elements
then their product $\lambda_i(\beta_0,\beta_0^\prime)\in\C$ is
given
by
$$
\lambda_i(\beta_0,\beta_0^\prime)
=(\beta_0,D_i\beta_0^\prime+D_{i-1}\beta_1^\prime+\dots+
D_1\beta_{i-1}^\prime)
$$
where $\beta_0^\prime, \beta_1^\prime, \dots,\beta_{i-
1}^\prime$
form a
solution of $(23_i)$.
\endproclaim

3.14. As an application of the above general results, consider
deformations
of flat bundles.

Suppose that $\nabla_t$, for $-\epsilon<t<\epsilon$, is an
analytic
family
of flat connections on a vector bundle $\E$ over a compact
Riemannian
manifold $M$ of odd dimension $2l-1$.
We assume that all $\nabla_t$ preserve a fixed
Hermitian metric on $\E$. Then
$$\nabla_t=\nabla +\sum_{i\ge 1}t^i\Omega_i,$$
where $\nabla=\nabla_0$ and
$\Omega_i\in A^1(M;\End(\E))$ are 1-forms with values in the
bundle
of endomorphisms. The latter bundle has a natural flat
structure
induced
by $\nabla$, and from $\nabla_t^2=0$ it follows that the first
form
$\Omega_1$ is flat $\nabla(\Omega_1)=0$. Thus, it determines a
cohomology
class
$$\omega_1=[\Omega_1]\in H^1(M;\End(\E)).$$
Consider now the following bilinear form
$$q:H^{l-1}(M;\E)\times H^{l-1}(M;\E)\to \C$$
given by
$$q(\alpha,\beta)=<\omega_1\cdot\alpha\cup\beta,[M]>,$$
where $\alpha,\beta\in H^{l-1}(M;\E)$, the product
$\omega_1\cdot\alpha$
belongs to $H^l(M;\E)$, and the cup-product $\cup$
uses the Hermitian metric on $\E$.

The form $q$ is $(-1)^l$-Hermitian; this fact can be easily
verified
using
the assumption that all $\nabla_t$ preserve the Hermitian
metric.
Let
$\sign(q)$ denote the signature of the form $q$, if $l$ is
even, and
the
signature of $iq=\sqrt{-1}\cdot q$, if $l$ is odd.

The following result essentially coincides with
the main result of \cite{\Kirk}.
\proclaim{3.15. Corollary} In the situation described above,
consider the
family of Atiyah-Patodi-Singer operators
$$B_t: A^{ev}(M;\E) \to A^{ev}(M;\E),\quad
-\epsilon<t<\epsilon,$$
where $B_t\phi=i^l(-1)^{p+1}(\ast \nabla_t -\nabla_t
\ast)\phi$
acting on
a form $\phi\in A^{2p}(M;\E)$. Then the first
signature $\sigma_1$ of the linking
form 3.8 associated with this analytic self-adjoint family of
elliptic
operators $B_t$ is equal to $\sign(q)$.
\endproclaim
\demo{Proof} Applying Corollary 3.11, we obtain that that the
first
signature
$\sigma_1$ can be computed as the signature of the pairing
$(x,y)\mapsto (B_1x,y)$, where $B_1=\pm(\ast\Omega_1-
\Omega_1\ast)$,
acting on the space
of harmonic forms of even degrees
with respect to the flat connection $\nabla$.
We can identify the above space of harmonic forms with
$H^{\ev}(M;\E)$. Our
aim now is to find a form on the middle dimensional cohomology
with
the
same signature.

We will use the following known lemma: {\it
let $l: X\times X\to\C$ be a Hermitian
form on a finite dimensional vector space $X$ and let
$A\subset B\subset X$ be subspaces such that $A^{\perp}=B+K$,
where
$K$ is the annihilator of $l$, $K=X^\perp$. Then the signature
of
the induced
form on $B/A$ is equal to $\sign(l)$.}

Consider first the case when the number $l$ is odd, i.e.
$l=2r-1$,
and so
the dimension of $M$ is $4r-3$. Let $A$ be the direct sum of
all
$H^{2k}(M;\E)$
with $k<r-1$; let $B$ denote the direct sum of all
$H^{2k}(M;\E)$
with
$k\le r-1$. Applying the above mentioned lemma to this
situation we
get that
the signature $\sigma_1$ is equal to the form on $H^{l-
1}(M;\E)$
given by
$$(\alpha,\beta)\mapsto
i\cdot(\ast\Omega\alpha,\beta)=i\int_M\Omega\alpha
\wedge\beta$$
as applied to the harmonic representatives. It is well-known
that
the last
pairing can be expressed in terms of the cup-product on
homology.
This proves our statement in the case of odd $l$.

In the case of $l$ even the arguments are similar. We apply
the
above lemma
to the subspace $A$, being the sum of all $H^{2k}(M;\E)$ with
$k>r$,
and
to the subspace $B$,
equal to the sum of all $H^{2k}(M;\E)$ with $k\ge r$. By the
lemma,
we
obtain that $\sigma_1$ is the signature of the form on
$H^{2r}(M;\E)$
given by
$$(\alpha,\beta)\mapsto
(\Omega\ast\alpha,\beta)=\int_M\Omega\ast\alpha
\wedge\ast\beta.$$
Identifying $H^{2r}(M;\E)$ with $H^{l-1}(M;\E)$ via $\ast$, we
arrive
to the statement of the Theorem.
$\square$
\enddemo

Higher terms of the spectral sequence of quadratic forms
determined
by the
family $B_t$ can also be computed using Corollary 3.13 in
terms of
some
Massey products; it can then be shown that these operations
are
correctly defined on the homology.
We will not continue along this way here, because our Theorem
1.5
gives a simple general answer in homological terms.

\heading 4. Parametrized Hodge decomposition \endheading

4.1. Suppose that we have a closed oriented Riemanian
manifold $M$ of dimension $n$
and a flat Hermitian vector bundle $\E$ of rank $m$ over $M$.
Assume that an analytic deformation of the flat structure of
$\E$ is
given;
this means
that we have an analytic one-parameter family of flat
connections
$$\nabla_t: A^k(M;\E) \to A^{k+1}(M;\E), \ \ k=0, 1,
\dots\tag24$$
where $t\in(-\epsilon,\epsilon), \ \nabla_t^2=0$ and the
Hermitian
metric on
$\E$ is flat with respect to every $\nabla_t$. Note that
analyticity
of the
family $\nabla_t$ we understand as in section 3.4; in the
present
situation
this means that $\nabla_t=\nabla_0+\sum_{i=1}^\infty
t^i\Omega_i$
where
$\Omega_i\in\A^1(M;\End(\E))$ and the power series converges
in any
Sobolev
space.

These data determine a germ of deformation of the twisted De
Rham
complex
$$\dots\to A^k(M;\E) @>{\nabla_t}>>A^{k+1}(M;\E)@>{\nabla_t}>>
A^{k+2}(M;\E) \to \dots$$
Instead we will study the following single cochain
complex of $\OO$-modules
$$\dots\to \OO A^k(M;\E) @>{\nabla}>>\OO
A^{k+1}(M;\E)@>{\nabla}>>
\OO A^{k+2}(M;\E) \to \dots\tag25$$
Here the symbol $\OO A^k(M;\E)$ denotes the set of germs of
analytic
curves
in $A^k(M;\E)$ (defined as in 3.8) and the map $\nabla$ acts
by the
formula
$$(\nabla\alpha)(t)=\nabla_t(\alpha(t))\tag26$$
for $\alpha\in \OO A^k(M;\E),\alpha:(-\epsilon,\epsilon)\to
A^k(M;\E),
t\in (-\epsilon,\epsilon)$.

The chain complex (25) will be called {\it the germ-complex of
the
deformation;}
is the central object of our study of
the deformation $\nabla_t$. The purpose of this section is to
prove
a version
of the Hodge decomposition theorem for this complex.

4.2. Recall first the operators which are defined on the
twisted DeRham complex. Every space $A^k(M;\E)$ of smooth $k$-
forms
carries a
Hermitian scalar product
determined by the metrics on $M$ and on $\E$. The Hodge
duality
operator
$$\ast: A^k(M;\E)\to A^{n-k}(M;\E)\tag27$$
satisfies
$$(\phi,\psi)=\int\limits_M \phi\wedge \ast\psi ,\ \
\phi,\psi\in A^k(M;\E),\tag28$$
where $(\phi,\psi)$ denotes the scalar product and $\wedge$
denotes
the exterior product
$$\wedge: A^k(M;\E)\times A^l(M;\E)\to A^{k+l}(M)\tag29$$
determined by the metric on $\E$ \ (here $A^{k+l}(M)$ denotes
the
space of
$(k+l)$-forms on $M$ with complex values). Note that this
exterior
product
satisfies
$$\alpha\wedge\beta =(-1)^{kl}\overline{\beta\wedge\alpha},\ \
\alpha\in
A^k(M;\E), \ \ \beta\in A^l(M;\E).$$

These structures on the twisted De Rham complex determine (by
pointwise
extension) the following objects on the germ-complex (25).
The scalar product (28) defines the following bilinear form
$$(\ ,\ ):\OO A^k(M;\E)\times \OO A^k(M;\E) \to \OO\tag30$$
where for germs of maps $\alpha,\beta:(-\epsilon,\epsilon)\to
A^k(M;\E)$
their product is given by
$$(\alpha,\beta)(t)=(\alpha(t),\beta(t)), \ \ t\in (-
\epsilon,\epsilon).$$
This form is $\OO$-linear in the first variable, Hermitian
$$(\alpha,\beta)=\overline{(\beta,\alpha)}$$
and positively defined; the latter means that the scalar
square
$(\alpha,\alpha)\in \OO$ is a germ of a {\it real\/} curve
assuming
{\it non-negative} values and $(\alpha,\alpha)=0$ if and only
if
$\alpha = 0$.

In a similar way the Hodge duality operator (27) and the
exterior
product (29) determine the following maps
$$\ast: \OO A^k(M;\E)\to \OO A^{n-k}(M;\E)\tag31$$
$$\wedge: \OO A^k(M;\E)\times \OO A^l(M;\E)\to \OO
A^{k+l}(M)\tag32$$
where
$(\ast\alpha)(t)=\ast(\alpha(t))$ for $\alpha\in\OO A^k(M;\E),
t\in
(-\epsilon,\epsilon)$ and
$(\alpha\wedge\beta)(t)=\alpha(t)\wedge\beta(t)$,
$\beta\in\OO A^l(M;\E)$.  The star-operator (31) is $\OO
A^0(M)$-
linear
and satisfies
$$\ast\ast\alpha=(-1)^{k(n-k)}\alpha$$
A relation between the star (31), the exterior product (32)
and the "scalar product" (30) is established by the formula
$$(\alpha,\beta)(t)=\int\limits_M (\alpha\wedge\ast\beta)(t)\
\
\tag33$$

Let
$$\nabla^{\prime}: \OO A^k(M;\E)\to \OO A^{k-1}(M;\E)\tag34$$
be given by
$$\nabla^{\prime}(\alpha)=(-
1)^{n(k+1)+1}\ast\nabla\ast(\alpha)$$
for $\alpha\in \OO A^k(M;\E)$. If $\beta\in\OO A^{k-1}(M;\E)$
then
$$(\alpha,\nabla\beta)=(\nabla^{\prime}\alpha,\beta)\ \in\OO$$
which means that $\nabla^\prime$ is dual to $\nabla$ with
respect to
the product (30).

Let
$$\Delta: \OO A^k(M;\E) \to \OO A^k(M;\E)\tag35$$
be the Laplacian $\Delta= \nabla\nabla^\prime +\nabla^\prime
\nabla$.

An element $\alpha\in \OO A^k(M;\E)$ is called {\it
harmonic\/} if
$\Delta\alpha=0$. Since
$$(\Delta\alpha, \alpha)=(\nabla\alpha,\nabla\alpha)+
(\nabla^\prime\alpha,\nabla^\prime\alpha)$$
and the scalar product is positively defined we obtain that a
form
$\alpha\in \OO A^k(M;\E)$ is harmonic if and only if
$\nabla\alpha=0$
and $\nabla^\prime\alpha=0$. The set of all harmonic forms in
$\OO
A^k(M;\E)$
will be denoted $\Har^k$. It is an $\OO$-module.

The following theorem is the main result of this section. In
its
statement
we use the notation introduced in 3.8: if $X\subset Y$ is a
$\OO$-submodule of an $\OO$-module $Y$ then
${\frak Cl}X$ denotes the set of all $y\in Y$ such that $fy$
belongs
to
$X$ for some nonzero $f\in \OO$.

\proclaim{4.3. Theorem} Suppose that an analytic deformation
$\nabla_t$ of
a flat
Hermitian vector bundle $\E$ over a compact manifold $M$
is given, cf. subsection
4.1. Consider the germ-complex of $\OO$-modules (25)
determined by
this
deformation. Then:

(1) the following decomposition holds:
$$\OO A^k(M;\E)
=\Har^k\ \oplus \ {\frak Cl}(\nabla(\OO A^{k-1}(M;\E)))
\ \oplus
\ {\frak Cl}(\nabla^\prime(\OO A^{k+1}(M;\E)))\tag36$$
and the terms of this decomposition are orthogonal to each
other
with respect to scalar product (30);

(2) the $\OO$-module $\Har^k$ of harmonic forms is free of
finite
rank;

(3) the factor-modules
$${\frak Cl}(\nabla(\OO A^{k-1}(M;\E)))/
\nabla(\OO A^{k-1}(M;\E))
\tag37$$
and
$${\frak Cl}(\nabla^\prime (\OO A^{k+1}(M;\E)))/
\nabla^\prime (\OO A^{k+1}(M;\E))
\tag38$$
are finitely generated torsion $\OO$-modules;

(4) the homology of complex (25) is finitely generated over
$\OO$
and is isomorphic to
$$\Har^k\ \oplus\ \tau^k$$
where
$$\tau^k=
{\frak Cl}(\nabla(\OO A^{k-1}(M;\E)))/
\nabla(\OO A^{k-1}(M;\E)) \tag39$$
is the torsion part of homology and $\Har^k$ is the free
part of the homology.

(5) the star operator (31) establishes $\OO$-isomorphisms
$$\Har^k\ @>{\simeq}>>\ \Har^{n-k},\ \ \
\tau^k\ @>{\simeq}>>\ \varrho^{n-k}\tag40$$
where we denote
$$\varrho^j=
{\frak Cl}(\nabla^\prime (\OO A^j(M;\E)))/
\nabla^\prime(\OO A^{j+1}(M;\E))\tag41$$

$\OO$-
that the
by
that

\endproclaim
Proof of a similar theorem in a more general situation of
elliptic complexes is given in \cite{\Fa}, \S 3.
The last statement (5) is absent in \cite{\Fa} (since it is
meaningless
for general elliptic complexes), but this statement obviously
follows
from the definitions.
$\square$

\heading 5. De Rham theorem for the germ-complex\endheading

In this section we are going to prove the
following statement which is similar to the classical De Rham
theorem.

\proclaim{5.1. Proposition} Let $\nabla_t$ be a deformation of
flat
Hermitian
structure and let the
$\OO [\pi]$-module $\VV$ be the deformation of the monodromy
representation as defined in subsection 1.3. Let
$H^{\ast}(M;\VV)$
be
the cohomology with coefficients in $\VV$, cf. subsection 1.4.
Then there is a canonical isomorphism between
$H^{\ast}(M;\VV)$
and the cohomology of the germ-complex (25). Moreover, if the
dimension
of the manifold $M$ is odd, $n=2l-1$, then the linking form
(6)
corresponds under this isomorphism to the $(-1)^l$-Hermitian
form
$$\tau^l \times \tau^l \to \M/\OO$$
which is given on the classes $[f],\ [f^\prime]\in \tau^l$ by
the
formula
$$\{[f],[f^\prime]\}(t)=\ t^{-k}\int\limits_M g(t)\wedge f^
\prime (t)\ \ \mod \OO\tag42$$
where $g\in\OO A^{k-1}(M;\E)$ is any solution of the equation
$t^kf\ =\nabla g$.
\endproclaim
The proof (cf. subsection 5.3 below) is based on the following
analytic
fact:

\proclaim{5.2. Lemma} Let $\E$ be a vector bundle over
a closed $n$-dimensional
ball $M$ lying in an Euclidean space $\R^n$ and let
$\nabla_t$, with
$t\in
(-\epsilon,\epsilon)$, be a family of flat connections on $\E$
which
is
analytic in sense of 3.4. Fix a point $p\in M$ and a vector
$e\in\E_p$ is
the fiber above $p$. For every $t\in (-\epsilon,\epsilon)$
there
exists a
unique section $s_t\in C^\infty(\E)$ such that
$\nabla_t(s_t)=0$ and
$s_t(p)=e$, cf. \cite{\Kob}, chapter 1. Then the curve of
sections
$(-\epsilon,\epsilon)\to C^\infty(\E)$, where $t\mapsto s_t$,
is
analytic
in the sense of 3.3.
\endproclaim
The proof given in 5.4.
\demo{5.3. Proof of Proposition 5.1} The arguments are
standard;
we will describe them briefly for completeness.

Let us define presheaves $\FF^i$, where $i=0,1,2,\dots$, on
$M$. For
an open
set $U\subset M$ let $\FF^i(U)$ denote
$\OO C^\infty((\Lambda^iT^\ast(M)\otimes \E)|_U)$;
in other words, the sections
of $\FF^i$ are germs of analytic curves of $i$-forms over $U$
with
values
in $\E$. Analyticity of curves of forms over an open set $U$
we
understand
in the following sense: for any compact
submanifold with boundary $C\subset U$,
the curves of restrictions to $C$ are supposed to be analytic
with
respect
to Definition 3.3. Note that $\FF^i(U)$ has natural structure
of a
$\OO$-module and the restriction maps of the presheaf are
$\OO$-
homomorphisms.

The path of flat connections $\nabla_t$ on $\E$ defines a
homomorphism
$\nabla:\FF^i\to\FF^{i+1}$ for every $i=0,1,2,\dots$, acting
as
follows:
let $\omega=(\omega_t)$ be an analytic curve of $i$-forms over
$U$;
then
$\nabla(\omega)$ is a curve of $(i+1)$-forms represented by
$t\mapsto \nabla_t(\omega_t)$. It is clear that $\nabla^2=0$.

Let $\BB$ denote the following presheaf:
$$\BB(U)=\ker[\nabla:\FF^0(U)\to\FF^1(U)].\tag43$$
According to Lemma 5.2, $\BB$ is a locally trivial presheaf of
flat
analytic curves
of sections of $\E$, i.e. analytic curves of sections $s_t$
such
that
$\nabla_t(s_t)=0$. In fact, if $U\subset M$ is a disk then
$\BB(U)$ is isomorphic to $\OO \E_p$ - the set of analytic
curves in
the fiber over a point $p\in U$.

We claim now that for any disk $U\subset M$ the sequence
$$0\to\BB(U)@>\nabla>>\FF^0(U)@>\nabla>>\FF^1(U)@>\nabla>>
\FF^2(U)@>\nabla>>\dots\tag44$$
is exact. To prove this, fix a point $p\in U$ and a frame in
the
fiber
$e_1, e_2,\dots, e_m\in\E_p$. By Lemma 5.2 we may find germs
of
curves
$s_1,\dots,s_m
\in \BB(U)$ such that $s_k(p)=e_k$ for $k=1,2,\dots, m$. Now,
any
$\omega\in\FF^i(U)$ can be uniquely represented as
$\omega=\sum_{k=1}^m \omega_ks_k$ where $\omega_k\in \OO
A^i(U)$. We
see that $\nabla(\omega)=0$ if and only if $d\omega_k=0$ for
all
$k=1,2,\dots,m$. The usual proof of the Poincare Lemma (cf.,
for
example,
\cite{\War}, 4.18) shows that if $\omega_k$ is an analytic
curve of
closed
differential $i$-forms on $U$ then there exists an analytic
curve
of $(i-1)$-forms $\nu_k$ with $d\nu_k=\omega_k$; the
construction of
such forms $\nu_k$ described in \cite{\War}, 4.18, uses
contraction
of the curve of forms $\omega_k$ with a
vector field and then an integration and so it produces curves
of
forms
analytic with respect to the parameter. Denote
$\nu=\sum_{k=1}^m
\nu_ks_k$.
Then $\nu\in\FF^{i-1}(U)$ and $\nabla \nu=\omega$.

Let $\UU=\{U_\alpha\}$ be a good finite cover of $M$, i.e. a
cover
by
interiors of closed disks such that all intersections of the
sets of
$\UU$ are also disks. Consider the \v Cech - De Rham complex
$K=\oplus K^{p,q}$, where $K^{p,q}=C^p(\UU,\FF^q)$
is the module of \v Cech $p$-cochains of the cover $\UU$
with coefficients in the presheaf $\FF^q$, cf. \cite{\Bo},
chapter
2.
The bigraded complex $K$ has differential $D=D^\prime
+D^{\prime\prime}$,
where $D^\prime$ is the \v Cech differential (it has degree
(1,0))
and
$D^{\prime\prime}$ acting on $K^{p,q}$ is equal $(-
1)^p\nabla$.
The differential $D^{\prime\prime}$  has degree (0,1).

Let $L_1$ denote the \v Cech complex $C^p(\UU,\BB)$ and let
$i_1:L_1\to K$
denote the natural inclusion $C^p(\UU,\BB)\to C^p(\UU,\FF^0)$;
let
$L_2$ denote the germ-complex (25) and let $i_2:L_2\to K$
denote the
natural inclusion $\OO A^q(M;\E)\to C^0(\UU;\FF^q)$. Since the
presheaves
$\FF^q$ are fine and the sequences (44) are exact over all
intersections of
the sets of $\UU$, the usual arguments, using the spectral
sequence
of the
bicomplex $K$, cf., for example \cite{\Bo}, chapter 3, give
that
both
maps $i_1$ and $i_2$ induce isomorphisms in cohomology. Thus,
the
cohomology
of the germ-complex is isomorphic to the sheaf cohomology
$H^\ast(\UU,\BB)$
which, as it is well known, is isomorphic to cohomology with
local
coefficients
$H^\ast(M,\VV)$, where $\VV$ is the local system determined by
the
monodromy
of the locally constant sheaf $\BB$. Obviously, $\VV$ is the
deformation
of the monodromy representation as defined in 1.3.

We are left to prove the second statement of Proposition 5.1,
concerning the
linking pairings. To do this we have to compare multiplicative
structures
on the germ-complex and on the \v Cech complex
$C^\ast(\UU,\BB)$.

Let $K^\prime$ denote the complex constructed similarly to $K$
with
respect to the trivial line bundle instead of $\E$ and the
exterior
derivative $d$ instead of the path of flat connections
$\nabla_t$.
In other words, $K^{\prime p,q}=C^p(\UU,\OO \Omega^q)$,
where $\Omega^q$ denotes the sheaf of differential forms.
Similarly, let $L_1^\prime$ denote $C^\ast(\UU,\OO)$
and $L_2^\prime$ denote the complex $\OO\Omega^\ast$ (the
germs of
curves
in the De Rham complex). We have the obvious imbeddings
$j_1:L_1^\prime\to K^\prime$ and
$j_2:L_2^\prime\to K^\prime$ which induce isomorphisms in
cohomology.

The multiplicative structure
in the \v Cech-De Rham complex $K$, which uses the Hermitian
scalar
product
on the bundle $\E$ to define the wedge product of form with
values
in $\E$,
and is given by a formula similar to \cite{\Bo}, p. 174,
determines
the
chain map
$$K\otimes \overline K\to K^\prime.\tag45$$
Here $\overline K$ denotes complex $K$ with the $\C$-module
structure
twisted by the complex conjugation.

Similarly, there are products $L_1\otimes \overline L_1\to
L_1^\prime$
and $L_2\otimes \overline L_2\to L_2^\prime$, which are
restrictions
of the
product (45). In other words, these products are related by
the
following commutative diagrams
$$
\CD
L_1\otimes\overline L_1@>>> L^\prime_1\\
@VV{i_1\otimes i_1}V @VV{j_1}V\\
K\otimes \overline K@>>> K^\prime,
\endCD
\qquad\quad
\CD
L_2\otimes\overline L_2@>>> L^\prime_2\\
@VV{i_2\otimes i_2}V @VV{j_2}V\\
K\otimes \overline K@>>> K^\prime
\endCD
$$

Let $\theta_2:L_2^\prime \to \OO$ be the $\OO$-homomorphism
given
by the fundamental class of $M$; in other
words
$$\theta_2(\omega_t)=\int_M \omega_t$$
(the above integral is defined to be zero if the degree of the
forms
$\omega_t$ is not equal to $n=\dim M$). $\theta_2$ is a
cocycle,
i.e. it
vanishes on coboundaries. Thus there exists a cocycle $\theta
:K^\prime
\to \OO$ of degree $n$ such that
$\theta|_{L_2^\prime}=\theta_2.$
Denote $\theta_1=\theta|_{L_1} : L_1^\prime \to \OO$;
it represents the fundamental class of $M$
in the \v Cech cohomology.

Each of these cocycles determines, in a standard way,
a linking pairing on the corresponding
complex. For example, the cocycle $\theta_1$ together with the
product
$L_1\otimes\overline L_1\to L^\prime_1$ determine the linking
pairing on the
\v Cech cohomology, which is identical with the pairing (7),
described
in subsection 1.4. The linking pairing corresponding to the
cocycle
$\theta_2$ and to the product $L_2\otimes\overline L_2\to
L^\prime_2$ is
the one which was desrcibed in the statement of Proposition
5.1.
They are
isomorphic to the linking pairing on the torsion subgroups of
the
cohomology
of $K$ which is constructed by using the product
$K\otimes\overline K\to K^\prime$ and the cocycle $\theta$.
$\square$
\enddemo

\demo{5.4. Proof of Lemma 5.2} We are thankful to V.Matsaev
who
suggested
the idea of the following arguments.

We may suppose that the vector bundle $\E$ has been
trivialized and
that the
connections $\nabla_t$ are of the form $d+\Omega$, where
$\Omega=\Omega(x,t)$ is a
matrix-valued 1-form on $M$ which depends on the parameter $t$
analytically
(in sense of definition 3.3). We will also suppose that $M$ is
the
ball
of radius 1 in $\R^n$ and that $p=0$ is its center.

Let $$Z=\sum_{i=1}^n r_i\frac{\partial}{\partial r_i}$$
be the radial vector field and
let the matrix valued function $F(x,t)$ be defined as
evaluation of
$\Omega(x,t)$ on the vector field $Z$. Then solving the system
of
linear differential equations of parallel displacement along
the
radius
joining a point $x\in M$ with the center $0\in M$,
we obtain the following formular
for the flat section $s_t(x)$
$$
\split
&s_t(x)=\\
&\sum_{k=0}^\infty (-
|x|)^k\int_{0\le\tau_k\le\dots\le\tau_1\le 1}
F(\tau_1x,t)\dots F(\tau_kx,t)
d\tau_1 \dots d\tau_k\cdot e
\endsplit
\tag46$$
This power series converges absolutely, cf. \cite{\Gan},
chapter 14,
\S 5.

Let $W\subset\C$ be a neighbourhood of the interval $(-
\epsilon,\epsilon)$
such that the curve of 1-forms $\Omega(x,t)$, where
$t\in (-\epsilon,\epsilon)$, can be extended to an analytic
curve
$t\mapsto A^1(M;\End(\E))$ defined for all $t\in W$, where
$A^1(M;\End(\E))$
is supplied with $L^2=\H_0$-topology.
Each term of the
series (46) is an analytic function mapping $W$ to
$\H_0(\End(\E))$;
thus
the series represents an analytic function and, therefore, the
curve
of
sections $s_t$ is analytic as a curve in $\H_0(\E)$.

We will show now that $s_t$ is analytic as a curve in any
Sobolev
space
$\H_k(\E)$. Let $\nabla$ denotes the connection $\nabla_t$ for
$t=0$.
Then the covariant derivative of $s_t$ along any vector field
$X$ is
equal to
$$\nabla_X(s_t)=-\Omega (X)\cdot s_t.\tag47$$
Sinse the right hand side is analytic as a curve in $\H_0(\E)$
we
obtain
that $s_t$ is analytic as a curve in $\H_1(\E)$.
For any other vector field $Y$ differentiationg (47), we
obtain
$$\nabla_Y\nabla_X(s_t)=-\nabla_Y(\Omega(X))\cdot s_t+
\Omega(X)\Omega(Y)\cdot s_t$$
By our assumptions the curve of 1-forms $\Omega$ and all its
derivatives
are analytic as curves in $\H_0$. Thus we obtain from the last
equation that
$s_t$ is analytic as curve in $\H_2(\E)$.

Continuing similarly, by induction we establish analyticity of
$s_t$
in the
sense of 3.3.$\square$
\enddemo

\heading 6. Proof of theorem 1.5\endheading

6.1. Suppose that $M$ is a closed oriented Riemannian manifold
of
odd
dimension $2l-1$ and $\E$ is a flat Hermitian vector bundle
over
$M$.
Suppose also that a deformation of the flat structure
of $\E$ is given. This means that we have a family of flat
connections
$\nabla_t$, preserving the
Hermitian metric on $\E$. Then the family of Atiyah-Patodi-
Singer
operators $B_t$ is defined
$$B_t\phi=i^l(-1)^{p+1}(\ast \nabla_t -\nabla_t \ast)\phi,
\ \ \ B_t: A^{ev}(M;\E) \to A^{ev}(M;\E),$$
all $B_t$ being elliptic and self-adjoint.
We want to compute the infinitesimal spectral flow of this
family.
According to our general results (cf. \S 3, Theorem 3.9)
we have to study the
signatures of the {\it analytic\/} linking form of this
deformation
which is given
by the general construction of $\S 3$ applied to the self-
adjoint
family $B_t$.
Our aim here is to compare this analytic linking form with the
{\it
algebraic\/} or {\it homological\/} linking form (6)
constructed in
1.4.

In general, the analytic linking form can be understood if we
know
the action
of the deformation on germs of holomorphic curves with values
in
$A^{ev}(M;\E)$ , cf (16).
The action of the Atiyah-Patodi-Singer operators $B_t$ on
holomorphic
curves with values in $A^k(M;\E)$ can be deduced from Theorem
4.3.
We will see that the information given by this theorem is
complete
enough for our purposes.

In order to simplify the notations we will denote the modules
appearing in Theorem 4.3 as follows
$$X^k=\nabla(\OO A^{k-1}(M;\E)),$$
$$Y^k=\nabla^\prime(\OO A^{k+1}(M;\E))$$
Then the decomposition of Theorem 4.3 looks
$$\OO A^k(M;\E)=\ \Har^k\ \oplus\ {\frak Cl}(X^k)\ \oplus\
{\frak
Cl}(Y^k)
\tag48$$
Recall the notation introduced in Theorem 4.3:
$${\frak Cl}(X^k)/X^k=\tau^k\qquad\text{and}\qquad
{\frak Cl}(Y^k)/Y^k=\varrho^k.$$
Note also that $\tau^k$ concides with the $\OO$-torsion part
of the
homology of the germ-complex (25) associated to the
deformation.

We find it convenient to consider separately the
cases of even and odd $l$.

6.2. Suppose first that $l$ is even, $l=2r$,
i.e. the dimension of $M$ is $4r-1$.

Consider the result of applying the family $B_t$ to
each term of decomposition (48). We obtain that the action
of the deformation $B_t$ on $\OO A^k(M;\E)$
splits into two sequences of {\it epimorphisms}:

$$(-1)^{r+p}\nabla\ast: \ {\frak Cl}(X^{2p})\ \to \ X^{4r-
2p}\tag49$$

$$(-1)^{r+p+1}\ast\nabla: \ {\frak Cl}(Y^{2p})\ \to \ Y^{4r-
2p-
2}\tag50$$

Applying the construction of subsection 3.8 to $B_t$
we obtain easily that the linking
form of this deformation splits into the following orthogonal
sum:
$$\bot_{p=0}^{r-1}(\tau^{2p} \oplus \tau^{4r-2p})\ \ \bot \ \
\bot_{p=0}^{r-1}(\varrho^{2p}\oplus \varrho^{4r-2p-2}) \ \bot
\
\tau^{2r}\tag51$$
where all forms $\tau^{2p} \oplus \tau^{4r-2p}$ and
$\varrho^{2p}\oplus \varrho^{4r-2p-2}$ for $0\le p \le {r-1}$
(except the last one on $\tau^{2r}$) are {\it hyperbolic} ,
cf.
subsection
2.12. Actually, their representation as sums of two Lagrangian
direct
summands is given by (51). Lemma 2.13 and the additivity of
the
signatures
$\sigma_i$ imply that
{\it the signatures of the analytic linking pairing
corresponding to
deformation of the Atiyah-Patodi-Singer operators $B_t$ are
equal to
the
signatures of the linking form
on the middle-dimensional torsion submodule $\tau^l$ given by
the
operator\/}
$\nabla\ast$, cf. (49).

Let us compute the last form explicitly.
By the construction of section 3.8 the
value of the form on holomorphic germs
$f,f^\prime\in\tau^{2r}$  is given by
$$\{f,f^\prime\}=t^{-k}(h,f^\prime)\ \mod\OO$$
where the germ $g$ is a solution of
the equation $t^kf=\nabla\ast h$. Using formulae (31) and (33)
we obtain that
$$\{f,f^\prime\}=t^{-k}\int\limits_M(\ast h\wedge f^\prime)\ \
\mod\OO$$
Using Theorem 5.1 we obtain from the last formula that {\it
the
analytic
linking form on $\tau^{2r}$ determined by the operator
$\nabla\ast$
coincides with the homological
linking pairing on the torsion part of homology as defined in
1.4.}

Now application of Theorem 3.9 finishes the proof of Theorem
1.5
in the case when $l$ is even.

6.3. In the case when $l$ is odd the arguments are similar.
Assume that $l=2r+1$; so the dimension of $M$ is $4r+1$.
As in the previous case we obtain that the
family of Atiyah-Patodi-Singer operators
$B_t$ acting on the germs of holomorphic curves, splits into
the
sequence of the following {\it $\OO$-epimorphisms\/}:
$$i(-1)^{r+p}\nabla\ast: \ {\frak Cl}(X^{2p})\ \to \ X^{4r-
2p+2},\tag52$$

$$i(-1)^{r+p+1}\ast\nabla: \ {\frak Cl}(Y^{2p})\ \to \ Y^{4r-
2p},\tag53$$
where $i=\sqrt{(-1)}$.
This implies that the linking form of self-adjoint family
$B_t$ is
the following orthogonal sum:
$$\bot_{p=0}^{r-1}(\tau^{2p} \oplus \tau^{4r-2p+2})\ \ \bot \
\
\bot_{p=0}^{r-1}(\varrho^{2p}\oplus \varrho^{4r-2p}) \ \bot \
\varrho^{2r}\tag54$$
Again, all forms in this decomposition
except the middle-dimensional form on $\varrho^{2r}$
(corresponding to the operator $(-i)\ast\nabla$) are {\it
hyperbolic} and so
have zero signatures by Lemma 2.13.
Thus, we obtain that the signatures of the linking
form corresponding to deformation of the Atiyah-Patodi-Singer
operator
are equal to the signatures of the linking pairing on
$\varrho^{2r}$
which acts as follows: for $f,f^\prime\in\varrho^{2r}$ their
product
$\{f,f^\prime\}\in\M/\OO$ is given by $\{f,f^\prime\}=t^{-
k}(g,f^\prime)$
where $g\in\OO A^{2r}(M;\E)$ solves the equation
$t^kf=-i\ast\nabla g$. Applying $\ast$ we get
$$t^k\ast f=-i\nabla g=\nabla(-ig)\tag55$$ and
$$\{f,f^\prime\}=i\times t^{-k}\int\limits_M (-ig)\wedge\ast
f^\prime=
i\times \{\ast f,\ast f^\prime\}^\prime\tag56$$
In (56) the brackets $\{\ ,\ \}^\prime$ denote the {\it
homological}
linking form
$$\{\ ,\ \}^\prime : \tau^{2r+1}\times\tau^{2r+1}\to\M/\OO$$
constructed as in subsection 1.4. The formula (56) shows that
the star-operator $\ast : \varrho^{2r}\to\tau^{2r+1}$
establishes an
isomorphism between the "analytic" form $\{\ ,\ \}$
on $\varrho^{2r}$ and the
algebraic form $\{\ ,\ \}^\prime$ on $\tau^{2r+1}$ multiplied
by
$i$.
The analytic form is Hermitian and the algebraic form is skew-
Hermitian;
taking into account our convention 2.15 on signatures of skew-
Hermitian
linking forms, we obtain finally that {\it the signatures
$\sigma_i$
of the
analytic linking form corresponding to the deformation of the
Atiyah-Patodi-Singer operator $B_t$, cf. $\S 3$,
coincide with the corresponding signatures of the algebraic
linking
form (7)
constructed in subsection 1.4.}

Application of Theorem 3.9 completes the proof of Theorem 1.5.
$\square$

\heading 7. Variation of the eta-invariant modulo $\Z$
\endheading

We now examine the behavior of the eta-invariant reduced
modulo
$\Z$.
This problem is rather well-understood, even (implicitly) in
the
original
work of Atiyah-Patodi-Singer \cite{\APS}; also see, for
example,
\cite{\CS}
and \cite{\Mathai}. Much more general results of this kind
were proven by P.Gilkey \cite{\Gilkey} (cf. Theorem 4.4.6 of
\cite{\Gilkey}
for example). For the convenience of the reader,
and in order to emphisize the explicit dependence on the
homotopy
type
of $M$, we give a complete and independent treatment here.
These
results,
together with Theorem 1.5, will be important for our study in
section 10 of
the problem of homotopy invariance of the $\rho$-invariant.

\proclaim{7.1. Theorem} Let $\E$ be a Hermitian line bundle
over a
closed
oriented Riemanian manifold $M$ of odd dimension $2l-1$.
Suppose
that two
flat connections $\nabla_0$ and $\nabla_1$ on $\E$ preserving
the
Hermitian metric are given. Let $\overline \eta_0$ and
$\overline
\eta_1$
denote the reductions modulo 1 of the eta-invariants of the
corresponding Atiyah--Patodi--Singer operators (2).
Consider the difference
$$\frac{1}{2\pi i}(\nabla_1-\nabla_0);$$
it is a closed 1-form on $M$ with real values.
Let $\xi\in H^1(M;\R)$ denote the corresponding cohomology
class.
Then the following formulae hold:
$$
\overline \eta_1-\overline \eta_0= \cases 0, &\text{if $l$ is
even}\\
2<\xi\cup L(M),[M]> \mod 1,&\text{if $l$ is odd},\endcases
$$
where $L(M)$ denotes the Hirzebruch polynomial in the
Pontrjagin
classes
of $M$.
\endproclaim

Note that the connections $\nabla_0$ and $\nabla_1$ are gauge-
equivalent
iff the class $\xi$ is integral, $\xi\in H^1(M;\Z)$. In this
case
we obviously obtain $\overline \eta_0 = \overline \eta_1 \mod
1$.

7.2. Theorem 7.1 can also be interpreted in the following way.
Assume that
$l$ is odd, $l=2r+1$, and so the dimension of $M$ is $4r+1$.
Let
$N_1, N_2,\dots, N_k$ be a set of oriented submanifolds of $M$
realizing
a basis in the homology group $H_{4r}(M;\Z)$; here $k$ is the
first
Betti number of $M$. Denote by $\tau_i$ the signature of
$N_i$,
$1\le i\le k$. Then Theorem 7.1 gives
$$\overline \eta_1-\overline \eta_0= 2\sum_{i=1}^k \tau_i x_i\
\mod
1$$
where the numbers $x_1, x_2,\dots,x_k$ are obtained as the
coefficients of decomposition
of $\xi\cap [M]$ (the Poincare dual of $\xi$) in term
of the basis formed by the classes $[N_i]\in H_{4r}(M;\Z):$
$$\xi\cap[M]=\sum_{i=1}^k x_i[N_i]\ \ \text{in} \ \
H_{4r}(M;\R)$$

7.3. Note that the space of flat structures on a given line
bundle
$\E$
up to gauge equivalence
is a torus $H^1(M;\R)/H^1(M;\Z)$; its dimension is equal to
the
first
Betti number of $M$. The eta-invariant $\overline \eta$
comprises a
function on this torus with values in $\R/\Z$. Theorem 7.1
states
that
this function is {\it linear\/} if $l$ is odd and is {\it
constant\/}
if $l$ is even. Theorem 7.1 implies also the following
statement:

\proclaim{7.4. Corollary} If $l$ is odd, $l=2r+1$, the reduced
eta-
invariant
$\overline\eta_\nabla \ \in \R/\Z$ is constant (i.e. does not
depend
on choice of the flat structure $\nabla$ on $\E$) if and only
if
all $4r$-dimensional compact submanifolds $N^{4r}\subset
M^{4k+1}$
have vanishing signatures.
\endproclaim
7.5. Consider now deformations of flat Hermitian bundles of
arbitrary
rank.

Again, let $M$ denote a compact oriented Riemannian manifold
of odd
dimension
$2l-1$ and $\E$ a vector bundle of rank $m$ over $M$. Let
$\det(\E)$
denote the line bundle $\wedge^m(\E)$, the $m$-th exterior
power of
$\E$.
Any connection $\nabla$ on $\E$ determines canonically a
connection
on
the line bundle $\det(\E)$ which will be denoted
$\det(\nabla)$.

Given a flat connection $\nabla$, the symbol
$\overline\eta_\nabla$
will denote the reduced modulo 1 eta-invariant of the
corresponding
Atiyah--Patodi--Singer operator (2).

\proclaim{7.6. Theorem} In the situation described above,
suppose
that we have two flat connections $\nabla_0$ and $\nabla_1$
which
can be
joined by a smooth path of flat connections $\nabla_t,\ 0\le
t\le 1$
on $\E$.
Then
$$\overline\eta_{\nabla_1}\ -\ \overline\eta_{\nabla_0}\ =\
\overline\eta_{\det(\nabla_1)}\ -\
\overline\eta_{\det(\nabla_0)}
\ \in\R/\Z$$
where $\det(\nabla_1)$ and $\det(\nabla_0)$ are the
corresponding
flat
connections on $\det(\E)$.

In particular, we obtain that if the dimension of $M$ is of
the form
$4r-1$, the eta-invariant $\overline\eta_\nabla$ assumes a
constant
value
(in $\R/\Z$) on connected components of the space of flat
connections.

If the dimension of $M$ is of the form $4r+1$ then
$$\overline\eta_{\nabla_1}\ -\ \overline\eta_{\nabla_0}\ =\
2<\xi\cup L(M),[M]>$$
where $\xi\in H^1(M;\R)$ is the cohomology class represented
by the
following closed 1-form
$$\frac{1}{2\pi i}(\det(\nabla_1)- \det(\nabla_0))=
\frac{1}{2\pi i}\tr(\nabla_1- \nabla_0);$$
thus, the reduced eta-invariant $\overline\eta_\nabla$ is
constant
on
connected components of the space of flat connections if and
only if
all $4r$-dimensional submanifolds $N^{4r}\subset M^{4r+1}$
have vanishing signatures.
\endproclaim

We refer to \cite{\Kob}, p.18  for general information
on determinant line bundles.

7.7. It can be useful in applications to express the class
$\xi$
which
appears in Theorems B and C in terms of the monodromy
representations.

Fix a base point $x\in M$. Let
$$\rho_\nu :\pi=\pi_1(M,x)\to U(\E_x),\ \ \nu=0, 1$$
denote the monodromy representation corresponding to the
connections
$\nabla_\nu$ where $\nu= 0, 1$. The monodromy representations
of the line bundle $\wedge^m\E_x$ corresponding to the
connection
$\det(\nabla_\nu), \ \nu=0,\ 1$,
is equal to the composition
$$\det\ \circ\ \rho_\nu:\ \pi \ \to\  U(\E_x)\ @>{\det}>>\
U(\wedge^m\E_x)$$
Let
$$\arg: U(\wedge ^m\E_x)\ \to \ \R/\Z$$
denotes the function argument. Then for any element $g\in \pi$
we
have
$$\frac{1}{2\pi }(\arg\ \det(\rho_1(g)) \ - \ \arg\
\det(\rho_0(g)))\ = \
-\ <\xi,g> \mod \Z.$$

The last equality determines the coset of class $\xi$ in
$H^1(M;\R)/H^1(M;\Z)$.

Results of Theorems 7.1 and 7.6 were obtained in
\cite{\Levine},
\cite{\Levinee}
in some special cases.

\subheading{7.8. Proofs of Theorems 7.1 and 7.6}
Suppose that we are in conditions of Theorem 7.6.
Namely, we suppose that $\E$ is a vector bundle of rank $m$
over
a closed oriented Riemannian manifold $M$ of odd dimension
$2l-1$
and $\nabla_t$ with $0\le t\le 1$ is a path of flat
connections on
$\E$.
Consider the vector bundle $\tilde\E$ over the product
$I\times M$
(where
$I$ denotes the interval $[0,1]$) induced from $\E$ by the
projection
$I\times M\to M$. The path of connections $\nabla_t$
determines a
unique
connection
$$\tilde\nabla : A^0(I\times M;\tilde \E)\to A^1(I\times
M;\tilde
\E)\tag57$$
on the bundle $\tilde\E$ where
$$(\tilde\nabla
s)(t,x)=(\nabla_ts(t,\cdot))(x)+dt\wedge\frac{\partial s(t,x)}
{\partial t}$$
for $s\in A^0(I\times M;\tilde \E)$, cf. \cite{\BGV}, p. 48.

Consider the {\it generalized signature operator\/} on the
manifold
$I\times M$ (supplied with the Riemannian metric which is the
product of the
metric of $M$ and the Euclidean metric on $I$)
$$A:\Omega^+\to\Omega^-$$
as defined by Atiyah, Bott and Patodi in \cite{\ABP}, p. 309.
Recall
that
here $\Omega^+\oplus\Omega^- =A^{\ast}(I\times M;\tilde E)$
and
$\Omega^{\pm}$ are the $\pm$-eigenspaces of the involution
$$\tau(\alpha)=i^{p(p-1)+l}\ast(\alpha),\ \
\alpha\in A^p(I\times M;\tilde \E)$$ and the operator $A$ is
$\tilde \nabla\ +\ \tilde\nabla^{\ast}$. Applying to this
operator
$A$
the index theorem of Atiyah--Patodi--Singer \cite{\APS}, part
I,
and taking
into account computation of the index density in \cite{\ABP},
section 6,
(cf. also \cite{\BGV}) we obtain the following equality:
$$\index A=2^l\int\limits_{I\times
M}\ch(\tilde\E)\cdot\L(I\times
M)[M]
-1/2(h_1-h_0)
-1/2(\eta_1-\eta_0).\tag58$$
Here $\ch(\tilde \E)$ denotes the Chern character form of the
connection
$\tilde\nabla$ on $\tilde \E$;\  $\L(M)$ denotes Hirzebruch
polynomial
in Pontrjagin forms which corresponds to
$\prod\frac{x_j/2}{\tanh(x_j/2)}$;
for $i=0, 1$ the numbers $h_i$ and $\eta_i$ denote
respectively
the dimension of the kernel and the eta-invariant of the
Atiyah--Patodi--Singer operator (2)
$$\pm(\ast\nabla_i-\nabla_i\ast)\ :A^\ast(M;\E)\to
A^\ast(M;\E),$$
acting of the full twisted De Rham complex.

Consider (58) modulo 1. Note that the numbers $h_i$ are equal
to the sum of all Betti numbers of $M$ with coefficients in
the flat
vector bundle $\E$ determined by the connection $\nabla_i$;
the alternating sum of those Betti numbers does not depend on
the
flat structure and is zero by Poincare duality. This proves
that {\it the numbers $h_i$ are even\/}, $i=0,\ 1$. Thus $h_i$
will
disappear from (58) if we consider it modulo 1.

Let $\eta_{\nabla_i}$ denote the eta-invariant of the
Atiyah--Patodi--Singer operator (2)
$$\pm(\ast\nabla_i -\nabla_i\ast): A^{ev}(M;\E)\to
A^{ev}(M;\E)\tag59$$
Then $2\eta_{\nabla_i}=\eta_i$ and we obtain from (58)
$$\eta_{\nabla_1}- \eta_{\nabla_0} =2^l\int\limits_{I\times M}
\ch(\tilde\E)\cdot\L(I\times M)\ \ \mod 1\tag60$$

Let us compute the Chern character $\ch(\tilde\E)$. The
curvature of
the connection $\tilde\nabla$ is equal to
$$K\ =\ dt\wedge\frac{d\nabla_t}{dt}.$$
Here $\frac{d\nabla_t}{dt}\in A^1(I\times M;\End(\E))$.
Thus we obtain
$$\ch(\tilde\E)=\ m \ -\ (2\pi i)^{-1}\tr(K)=
\ m\ -\ (2\pi i)^{-1}dt\wedge\tr(\frac{d\nabla_t}{dt})$$
Substituting this into (60) and using the fact that
$$\int\limits_{I\times M} \L(I\times M)= 0$$
(since the Pontrjagin forms of $I\times M$ do not contain
$dt$) we
obtain
$$\eta_{\nabla_1}- \eta_{\nabla_0} =\ - \frac{2^{l-1}}{\pi i}
\int\limits_{I\times M} dt\wedge
\tr(\frac{d\nabla_t}{dt})\wedge
\L(I\times M)\ \ \mod 1\tag61$$
Integrating the last formula with respect to $t$ we get
$$\eta_{\nabla_1}- \eta_{\nabla_0} =\ \frac{2^{l-1}}{\pi i}
\int\limits_M\tr(\nabla_1-\nabla_0)\wedge \L(M)\ \ \mod
1\tag62$$
Note that here $tr(\nabla_1-\nabla_0)\in A^1(M)$. This shows
that
only
the component of $\L(M)$ having dimension $2l-2$ appears in
(62).
Thus the LHS of (62) vanishes for $l$ even.

If $l$ is odd, $l=2r+1$, then we have
$$2^{2r}\L(M)_{4r}=\ L(M)_{4r},\tag63$$
where $L(M)$ denotes the Hirzebruch polynomial in the
Pontrjagin
classes defined by the generating series $x/\tanh(x)$ and the
subscript
refers to the corresponding homogeneous components.
On the other hand,
$$\det(\nabla_1)-\det(\nabla_0)=\tr(\nabla_1-\nabla_0)\tag64$$
(cf. \cite{\Kob}, p.18);
thus, $\tr(\nabla_1-\nabla_0)$ is a closed 1-form realizing
class
$2\pi i\xi$, cf. Theorem 7.1. This gives
$$\eta_{\nabla_1}- \eta_{\nabla_0}=2<\xi\cup L(M),[M]>\ \mod
1$$
and finishes the proof of Theorems 7.1 and 7.6.
$\square$

\heading 8. An example\endheading
In this section we consider the simplest possible example of
the
circle. This example was calculated analytically in
\cite{\APS}, II,
pages 410-411.
We wish to apply our theorems 1.5 and 7.1 in order to
illustrate
them.

Let $M$ be $S^1$, the circle, and let $\E$ be the trivial line
bundle over
$M$. Suppose that $\nabla_a$ is an analytic family of flat
connections
on $\E$ defined for $0\le a<1$ such that the induced family of
the
monodromy representations is given by
$$\rho_a:\pi\to U(1)=S^1,\ \ \rho(\tau)=\exp(2\pi i a)$$
where $\pi=\pi_1(M)$ and $\tau\in\pi$ is a generator.  If
$\eta_a$
denotes
the eta-invariant of the corresponding Atiyah--Patodi--Singer
operator,
which in the present case is
$$-i\ast\nabla_a : C^\infty \to C^\infty,$$
then the computation in \cite{\APS} gives:
$$
\eta_a=\cases 0,&\text{if $a=0$}\\
1-2a,&\text{if $0 < a < 1$}. \endcases
$$
Thus, for this family of flat connections the corresponding
monodromy
representations are parametrized by the circle and the eta-
invariant
has a jump near the trivial representation (which corresponds
to the
value
$a=0$) and it is smooth near all other representations.
Note that, near the trivial representation we have
$$\eta_0=0,\ \ \eta_+=1,\ \ \eta_-=-1$$
where the notation introduced in section 1 is used.

Let us compute the value of the jump using Theorem 1.5. In
order to
do this
we have to:

($i$) find the corresponding deformation of the monodromy
representation;

($ii$) find the cohomology with local coefficients
$H^{\ast}(\Hom
_{\Z[\pi]}
(C_{\ast}(\tilde M),\VV))$;

($iii$) calculate the linking form (6) on the middle-
dimensional
torsion;

($iv$) find the signatures $\sigma_i, \ i\ge 1$, cf. section
2.

The deformation of the monodromy representation near a
representation
$\rho_a$ is given by $\OO[\pi]$-module $\OO$ with the action
of
$\pi$
determined by
$$\tau \cdot f=\exp(2\pi i (t+a))\cdot f,\ \
t\in(-\epsilon,\epsilon),\ \
f\in \OO$$
where $\tau \in\pi$ is a generator.

To find the cohomology with local coefficients ($ii$) consider
the
cell
decomposition of $M=S^1$ consisting of one zero-dimensional
cell
$e^0$ and
one one-dimensional cell $e^1$. Then the chain complex
$C_\ast(\tilde M)$
of the universal cover $\tilde M$ is
$$0\to \Z[\tau,\tau^{-1}]e^1@>d>> \Z[\tau,\tau^{-1}]e^0\to 0$$
where
$d(pe^1)=(\tau-1)pe^0$ for $p\in \Z[\tau,\tau^{-1}].$
Thus, the cochain complex
$$\Hom
_{\Z[\pi]}
(C_{\ast}(\tilde M),\VV)$$
in this case is
$$0@<<< \OO@<\delta<<\OO@<<<0,$$ where
$\delta(f)=\exp(2\pi i(t+a))-1)\cdot f$ for $f\in\OO$.
We obtain that for $a\ne 0$ the last complex is acyclic,
confirming
that
{\it there are no jumps near nontrivial representations}.

In the case $a=0$ we obtain
$$H^1=\C=\OO/t\OO;$$
if $\alpha$ denotes the generator of the last group then the
linking
form
is given by
$$\{\alpha,\alpha\}=(\exp(2\pi i t)-1)^{-1}\in\M/\OO.$$
The linking form is skew-Hermitian (in this dimension) and in
order
to
compute its signatures we have first to multiply it by
$i=\sqrt{(-
1)}$,
cf. 2.15. Then
(using the notation of section 2) we  obtain
$$[\alpha,\alpha]=\Res i\{\alpha,\alpha\}=(2\pi)^{-1}$$
and thus
$$\sigma_1=1, \ \ \text{and}\ \ \sigma_j=0\ \ \text{for}\ \
j>1$$
Theorem 1.5 gives now the correct jump formulae near the
trivial
representation.

Theorem 7.1 together with remark 7.7 give the correct
reduction
of the eta-invariant modulo 1.

\heading 9. Deformations of flat line bundles and Blanchfield
pairings
\endheading

In this section we present a more general example of
application of
our
Theorem 1.5. We will show that the signatures derived
from studying the Blanchfield pairings, which are among the
most
standard
tools of the knot theory, are special cases of the signatures
studied in \S 3, corresponding to some particular curves of
deformations of line bundles.

We will first recall the construction of the
Blanchfield pairing (introduced by Blanchfield [B]) and its
local
version in the form convenient for the sequel. The result of
this
section
follows from the work of W.Neumann \cite{\Ne}
who studied the eta-invariant in this particular situation.

9.1. Let $M$ denote a compact oriented manifold of odd
dimension $2l-1$ and let
$$\phi: \pi_1(M)\ \to \ \Z\tag65$$
be a fixed epimorphism. Here $\Z$ will be understood as
the multiplicatively written infinite
cyclic group, whose generator will be denoted by $\tau$.
We will also fix a subfield $K$ in the field of complex
numbers.

Consider the infinite cyclic covering $\tilde M$ of $M$
corresponding
to the kernel of the homomorphism $\phi$.
Let $C$ denote the simplicial chain complex (with coefficients
in
$K$)
of $\tilde M$. $C$ is a complex of free finitely generated
left $\Lambda$-modules, where $\Lambda=K[\Z]=K[\tau, \tau^{-
1}]$
is the ring of
Laurent polynomials of $\tau$. For $0\le k\le 2l-1$, the
cohomology
$$H^k(M;\Lambda)=H^k(\Hom_{\Lambda}(C,\Lambda))\tag66$$
is a finitely generated $\Lambda$-module and (since $\Lambda$
is a
principal
ideal domain) we have the following decomposition
$$H^k(M;\Lambda)=\T^k\ \oplus\ F^k\tag67$$
where $\T^k$ denotes the torsion part of the $\Lambda$-module
$H^k(M;\Lambda)$ and $F^k$ denotes its free part.

Denote by $\RR$ the field of rational functions of $\tau$ with
coefficients
in $K$. We will consider $\RR$ with the involution
$\RR \to\RR$ which is the composition of the complex
conjugation and
substitution $\tau\mapsto\tau^{-1}$. We will denote this
involution
by
the overline. The ring $\Lambda$ is embedded into
$\RR$ and the involution preserves $\Lambda$.

The {\it Blanchfield form} is the map
$$\{\ ,\ \}: \T^l\ \otimes\ \T^l\ \to \ \RR/\Lambda\tag68$$
which is constructed as follows.
The extension $$0\to\Lambda\to\RR\to\RR/\Lambda\to 0$$
generates the
exact sequence
$$\dots \to H^{l-1}(M;\RR)\to H^{l-
1}(M;\RR/\Lambda)@>{\delta}>>H^l(M;\Lambda)
\to H^l(M;\RR)\to\dots$$
where $\delta$ is the Bockstein homomorphism. The image of
$\delta$
is
precisely the torsion submodule $\T^l$. For $\alpha,\beta\in
\T^l$
one
defines
$$\{\alpha,\beta\}=<\delta^{-1}(\alpha)\cup\beta, [M]>\
\in \RR/\Lambda\tag69$$
where the cup-product is taken with respect to the coefficient
pairing
$$\RR/\Lambda\times\Lambda\to \RR/\Lambda, \ (f,g)\mapsto
f\cdot
\overline g,\tag70$$
the dot stands for the multiplicaion of functions.

It is well known that the resulting pairing is {\it well-
defined,
non-degenerate\/} and $(-1)^l$-{\it Hermitian\/}
with respect to the induced involution on
$\RR/\Lambda$.

9.2. The {\it local version\/} of the Blanchfield pairing (68)
corresponding to a
is a prime ideal $\p$ in $\Lambda$,
is given by the map
$$\{\ ,\ \}_{\p}: \T^l_{\p}\ \otimes\ \T^l_{\p}\ \to \
\RR/\Lambda_{\p}
\tag71$$
Here $\Lambda_{\p}$ denotes the localization of the ring
$\Lambda$
with respect to the complement of the ideal $\p$ and
$\T^l_{\p}$
denotes
the $\p$-torsion part of the cohomology module
$H^l(M,\Lambda)$ i.e. the set of elements $z\in
H^l(M,\Lambda)$ such
that
for any $q\in \p$ there exists $n$ such that $q^nz=0$.
Clearly,
$\T^l_{\p}$ is a module over $\Lambda_{\p}$.

The prime ideal $\p$ is supposed to be
invariant under the involution.
The local Blanchfield pairing (71) is defined as the
restriction of
the
pairing (68) on the $\p$-torsion subgroup $\T^l_{\p}$ and then
reducing the
values modulo $\Lambda_\p$, i.e.
$$\{\alpha,\beta\}_\p=\{\alpha,\beta\} \mod \Lambda_\p $$

The global Blanchfield pairing (68) is direct orthogonal sum
of
local
pairings (71).

9.3. In the case when the field $K$ is $\C$, the field of
complex
numbers,
the prime ideals $\p$ in $\Lambda$ are in one-to-one
correspondence
with
complex numbers $\xi\in\C,\ \xi\ne 0$; the point $\xi\in\C$
represents
the principal ideal $\p\subset\Lambda$ generated by $\tau-
\xi$.

In this case we will write $\T^l_{\xi}$ instead of $\T^l_{\p}$
and
$\Lambda_\xi$ insead of $\Lambda_\p$. Thus the local
Blanchfield
form
in case $K=\C$ is denoted by
$$\{\ ,\ \}_{\xi}: \T^l_{\xi}\ \otimes\ \T^l_{\xi}\ \to \
\RR/\Lambda_{\xi}
\tag72$$

9.4. Now we are going to apply apply Theorem 1.5 in the
following
situation.
Suppose that $M$ is a compact manifold of odd dimension $2l-1$
and
$\phi:\pi_1(M)\to \Z$ is a fixed epimorphism.
Consider the following loop of one-dimensional
representations:
$$\rho_t: \pi\ =\ \pi_1(M)\to S^1=U(1),\qquad \rho_t\
=\ \mu_t\circ\phi\tag7.3$$
where
$$\mu_t:\Z\to  S^1,\quad \mu_t(\tau)=\exp(2\pi i t), \quad
0\le t\le
1.
\tag74$$
Here as above $\Z$ denotes the infinite cyclic group written
multiplicatively
and $\tau$ denotes its fixed generator.

Let $\omega$ be a closed 1-form on $M$ with real values
representing
the
De Rham cohomology class determined by $\phi$. In other word,
$\omega$
has the following property: for any closed loop $\alpha$ in
$M$
$$\phi([\alpha])=\tau^l\quad\text{where}\quad l=
\int\limits_\alpha
\omega.$$
Then for every $0\le t\le 1$
the operator
$$\nabla_t\ =\ d\ -\ 2\pi it\omega\wedge$$
is a flat connection on the trivial complex line bundle $\E$
over
$M$.
Note that the monodromy of the connection $\nabla_t$ is
$\rho_t$.

Thus, we have an analytic curve of flat connections.
We intend to compute the jumps of the eta-invariant
$\eta_t$ by using Theorem 1.5.

9.5. Fix a point $\xi\in S^1\subset \C$ on the unit circle.
It determines the
following representation
$$\nu_\xi:\pi\ \to\
U(1),\quad\nu_\xi(g)=\xi^n\quad\text{where}
\quad g\in\pi\ \text{and}\ \phi(g)\ =\ \tau^n\tag75$$

Then the 1-parameter family of representations
$$\mu_t\ =\ \nu_{\xi\exp(2\pi i t)}\ :\ \pi_1(M)\ \to\ U(1),
\qquad -\epsilon<t<\epsilon\tag76$$
is a deformation of the representation $\nu_\xi$.

As explained in (1.3), the deformation of representation
determines
a
module over the group ring $\OO[\pi]$ with
coefficients in $\OO$. Since all representations under
consideration
factor through
$\phi:\pi\to\Z$, it is enough for our purposes to consider the
following
$\OO[\Z]=\OO[\tau,\tau^{-1}]$-module $V_\xi$.
Here $V_\xi$ is equal to
$\OO$ as an $\OO$-module and the action of $\tau$ is given by
the
formula
$$(\tau\cdot f)(t)\ =\ \xi\exp(2\pi
it)f(t)\qquad\text{for}\quad
f\in\OO.$$

Then the cohomology of $M$ with coefficients in $V_\xi$ can be
computed
as cohomology of the following complex
$$\Hom
_{\OO[\tau,\tau^{-1}]}
(C(\tilde M); V_\xi)$$
where $\tilde M$ denotes the space of the infinite cyclic
covering
corresponding to the kernel of the homomorphism $\phi$ and
$C(\tilde M)$ denotes the simplicial chain complex of $\tilde
M$.

We want to show that {\it the linking form of this particular
deformation $\mu_t$ essentially
coincides with the local Blanchfield form\/} (72). More
precisely
this result
is formulated in 9.8. Note that W.Neumann \cite{\Ne} already
shown
that
the Blanchfield pairing determines the $\eta$-invariant in
this
case.

Note, that we have the following ring homomorphism
$$\alpha_\xi:\Lambda\to\OO,\qquad\text{where}\quad
\tau\mapsto\xi\exp(2\pi it)\in\OO,\tag77$$
and $V_\xi$ is just $\OO$
considered as a $\Lambda$-module via $\alpha_\xi$.  We observe
next
that $\alpha_\xi$ has a unique extension to
$\alpha_\xi:\Lambda_\xi\to\OO$
and also $\alpha_\xi:\RR\to\M$. The last map is a field
extension.

We have the following commutative diagram
$$
\CD
\Lambda@>>>\Lambda_\xi@>>>\RR@>>>\RR/\Lambda_\xi\\
       @. @V{\alpha_\xi}VV@V{\alpha_\xi}VV @VV{\alpha_\xi}V\\
       @. \OO@>>>\M @>>>\M /\OO
\endCD
\tag78$$

\proclaim{9.6. Lemma} $V_\xi$ is flat as a $\Lambda$-module.
\endproclaim
\demo{Proof} Since $\Lambda$ is a principal ideal domain it is
enough to
show that $V_\xi$ has no $\Lambda$-torsion. But this can be
easily
checked.
$\square$
\enddemo

\proclaim{9.7. Lemma} Let $X$ be a finitely generated
$\Lambda$-
module and
let $X_\xi$ denote its $\xi$-torsion, i.e.
$X_\xi=\ker((\tau-\xi)^n:X\to X)$, where $n$ is large.
Then the map
$$X_\xi\to X\otimes_{\Lambda} V_\xi,\qquad x\mapsto x\otimes
1\tag79$$
establishes an isomorphism between $X_\xi$ and the $\OO$-
torsion
submodule
of $X\otimes_{\Lambda} V_\xi$.
\endproclaim
\demo{Proof} Since the statement is clear when $X$ is free, it
is
enough
to prove it in the case when $X=\Lambda/(\tau-\xi)^n\Lambda$.
Then
$X\otimes V_\xi$ is isomorphic to
$$\OO/(\exp(2\pi it)-1)^n\OO\ \simeq\ \OO/t^n\OO$$
and the map
$$\Lambda/(\tau-\xi)^n\Lambda\ \to \ \OO/t^n\OO,\qquad
\tau\mapsto
\xi\exp(2\pi it)$$
is an isomorphism.  $\square$
\enddemo

\proclaim{9.8. Proposition} The linking form of the
deformation
$\nu_t$
coincides with the composition
$$\T^l_\xi\times\T^l_\xi\ \to\ \RR/\Lambda_\xi
\overset\alpha_\xi\to
\longrightarrow\M/\OO,\tag80$$
where $\T^l_\xi$ is the $\xi$-torsion submodule of the
$\Lambda$-module $H^l(M,\Lambda)$ , the first map is the
Blanchfield
pairing localized at $\xi$\ (cf. (72)), and the map
$\alpha_\xi$ is
given
by $f\mapsto f(\xi\exp(2\pi i t)),\quad f\in\RR$.
\endproclaim
\demo{Proof} Using Lemma 9.6 we obtain the following
isomorphisms
$$H^l(M;V_\xi)\simeq H^l(\Hom(C,V_\xi))\simeq
H^l(\Hom(C,\Lambda)\otimes
V_\xi)\simeq H^l(M;\Lambda)\otimes V_\xi$$
where $C$ denotes the simplicial chain complex of $\tilde M$
and
the tesor product is taken over $\Lambda$.
Then using Lemma 9.7 we may identify the $\OO$-torsion in
$H^l(M;V_\xi)$
with the $\xi$-torsion in $H^l(M;\Lambda)$, i.e. with
$\T^l_\xi$.

Two Bockstein homomorphisms (one, which is used in the
definition of
the
Blanchfield form and the other which is used in the
construction of
the
linking form of the deformation) appear in the following
commutative diagram
$$
\CD
H^{l-1}(M;\RR/\Lambda_\xi)@>{\delta}>>H^l(M;\Lambda_\xi)\\
@V{\alpha_\xi}VV @VV{\alpha_\xi}V\\
H^{l-1}(M;\M/\OO)@>{\delta}>>H^l(M;V_\xi)
\endCD
$$
There is also a commutative diagram
$$
\CD
\RR/\Lambda_\xi\times\Lambda_\xi@>>>\RR/\Lambda_\xi\\
@V{\alpha_\xi\times\alpha_\xi}VV@VV{\alpha_\xi}V\\
\M/\OO\times\OO@>>>\M/\OO
\endCD
$$
where the horizontal maps are the pairings which are used in
the
constructions of the cup-products used in the definition of
the
local Blanchfield form and the linking form of the
deformation. All
these facts taken together complete the proof.
$\square$
\enddemo

\heading 10. On homotopy-invariance of the $\rho$-invariant
\endheading

We now apply our results to the $\rho$-invariant.

10.1. Recall the definition
from \cite{\APS}. If $M$ is a closed smooth oriented odd-
dimensional
manifold and $\alpha$ a $k$-dimensional unitary representation
of
$\pi_1(M)$, then
$$\rho_{\alpha}(M)\ =\ \eta_{\alpha}(M)\ -\ \eta_0(M)$$
where $\eta_{\alpha}(M)$ denotes the eta-invariant of a flat
connection
with monodromy $\alpha$, and 0 denotes the trivial $k$-
dimensional
representation. The definition requires a choice of Riemannian
metric
on $M$ but it is shown in \cite{\APS}, as an easy consequence
of the
Index Theorem, that $\rho_{\alpha}(M)$ is independent of this
choice.
\cite{\APS} poses the problem of finding a direct
"topological"
definition
of $\rho_{\alpha}(M)$. In the special case that $(M,\alpha)$
bounds,
i.e. $M\ =\ \partial V$, where $V$ compact and oriented and
$\alpha$
extends
to a unitary representation $\beta$ of $\pi_1(V)$, then
$$\rho_{\alpha}(M)\ =\ k\sign(V)\ - \ \sign_{\beta}(V)$$
where $\sign_{\beta}(V)$ is the signature of the intersection
pairing on the
homology of $V$ with twisted coefficients defined by $\beta$,
and
$\sign(V)$
is the usual signature of $V$ -- see \cite{\APS} for more
details.

The related question of when $\rho_{\alpha}(M)$ is an oriented
homotopy
invariant of $M$ has been of some interest. W.Neumann
\cite{\Ne}
showed this
to be true when $\alpha$ factors through a free abelian group
and
S.Weinberger \cite{\We} extended this to a much larger class
of
torsion-free
groups. However it was already shown in \cite{\Wa} that the
$\rho$-
invariant
could distinguish (even mod $\Z$) homopoty equivalent fake
lens
spaces.
S.Weinberger
\cite{\We} finds such examples for a larger class of
fundamental
groups with
torsion.

Our goal is to apply Theorems 1.5 and 7.6 to give an explicit
homotopy invariant
definition of $\rho_{\alpha}(M)$, but with an indeterminacy
expressed by
a function on the space of unitary representations of
$\pi_1(M)$,
which is
constant on connected components -- zero at the component of
the
trivial
representation. If the representation space is connected (e.g.
$\pi_1(M)$
is free or free abelian) then the definition is complete. The
homotopy
invariance of $\rho_{\alpha}(M)$ up to a similar, but
rational-
valued,
indeterminacy was proved in \cite{\We} when $\pi_1(M)$
satisfies the
Novikov conjecture -- we will also recapture this result.

10.2. We will treat the $\rho$-invariant for a fixed closed
oriented
smooth
manifold $M$ of odd dimension as a function
$$\rho(M):\ \RR_k(\pi)\to\R,\qquad \rho(M)\cdot \alpha\ =\
\rho_{\alpha}(M),$$
where $\RR_k(\pi)$ is the real algebraic variety of $k$-
dimensional
unitary
representations of $\pi\ =\ \pi_1(M)$.
This function is piecewise continuous in the sense that there
is a
stratification of $\RR_k(\pi)$ by subvarieties
$$\RR_k(\pi)\ =\ V_0\supseteq V_1\supseteq V_2\dots$$
such that $\rho(M)|V_i-V_{i+1}$ is continuous (see
\cite{\Levinee}).
In fact, the discontinuities are integral jumps so that the
reduced
function
$$\overline\rho(M):\ \RR_k(\pi)\to\R/\Z$$
is continuous. We can rephrase Theorem 7.6 to give an explicit
description
of $\overline\rho(M)$ up to an indeterminacy expressed as a
function
on
$\RR_k(\pi)$ which is constant on connected components and 0
at the
trivial
representation. We will henceforth refer to such functions
(with
values in
$\R$ or $\R/\Z$) as {\it quasi-null}.

10.3. Define
$$\tilde\rho(M):\ \RR_k(\pi)\to\R/\Z$$
by the formula:
$$
\tilde\rho(M)\cdot\alpha\ =\ \cases
-2<(\arg \det \alpha)\cup\tilde L^{4r}(M),[M]>&\text{if} \dim
M=4r+1\\
0&\text{if} \dim M \equiv 3\mod 4.\endcases
$$
We explain the terms in this formula: $\arg\det \alpha$ is the
element
of $H^1(M;\R/\Z)$ defined by the homomorphism $\pi\to S^1$
given by
$g\mapsto \det \alpha(g),\ g\in\pi$, and $S^1\approx \R/\Z$
given by
$\exp(2\pi i t)\leftrightarrow t$. $\tilde L^{4r}(M)$ is a
lift to
$H^{4r}(M;\Z)$ of $L^{4r}(M)$, the Hirzebruch polynomial in
the
Pontjagin
classes of $M$. The existence of the integral lift was proven
by
Novikov
\cite{\No}. More explicitly, he observed that for any $\xi\in
H^1(M;\Z)$
$$<L^{4r}(M)\cup\xi,[M]>\ =\ \sign(N)$$
where $N$ is any closed oriented submanifold of $M$ dual to
$\xi$
(and so
$\dim N= 4r$). Note that $\tilde L^{4r}(M)$ is not unique, but
can
be varied
by a torsion class in $H^{4r}(M;\Z)$. This may change
$\tilde\rho(M)\cdot
\alpha$ by an element of finite order, which depends
continuously on
$\alpha$. Since $\tilde\rho(M)\cdot 0=0$ unambiguously, we
conclude
that
$\tilde\rho(M)$ is well defined up to a quasi-null function
with
values
in $\Q/\Z$.

We also point out that Novikov's formula gives an alternative
definition of
$\tilde\rho(M)$. Choose classes $z_1,\dots,z_n\in H_1(M;\Z)$
which
define
a basdis of $H_1(M;\Z)/torsion$ and let $z_1^\prime,\dots,
z_n^\prime$
be the dual basis of $H^1(M;\Z)$. Choose $N_1,\dots,N_n$
closed
oriented
submanifolds of $M$ Poincare dual to $z_1^\prime,\dots,
z_n^\prime$.
Then we have
$$\tilde\rho(M)\cdot\alpha\ =\ -2\sum_{i=1}^n
\sign(N_i)\arg\det
\alpha(z_i)$$
The choice of $z_1,\dots,z_n$ produces the indeterminacy of
$\tilde\rho(M)$ from this point of view.

Finally we point out that $\tilde\rho(M)$ is an oriented
homotopy
invariant
of $M$. This follows from the result of Novikov \cite{\No}
that for
any
$\xi \in H^1(\pi;\R)$, the invariant
$$(M,\phi)\ \mapsto\ <\phi^{\ast}(\xi)\cup L^{4r}(M),[M]>,$$
where $\phi :\pi_1(M)\to \pi$, is a homotopy invariant of
$(M,\phi)$.
If we choose an integer $m$ so that
$mH^1(M;\R/\Z)\subseteq H^1(M;\R)/H^1(M;\Z)$, then
$m\arg\det\alpha$
can be lifted to $H^1(M;\R)$ for any $\alpha$ and Novikov's
result
implies
that $m\tilde\rho(M)$ is a homotopy invariant of $M$. Thus, if
$M^\prime$
is homotopy equivalent to $M$, then $m(\tilde\rho(M)-
\tilde\rho(M^\prime))=0$
and so $\tilde\rho(M)-\tilde\rho(M^\prime)$ is a continuous
function
into
a discrete subset of $\R/\Z$. Thus $\tilde\rho(M)-
\tilde\rho(M^\prime)$ is a
quasi-null function into $\Q/\Z$.

We now rephrase Theorem 7.6 as:
\proclaim{10.4. Theorem} The difference
$$\tilde\rho(M)\ -\ \overline\rho(M):\ \RR_k(\pi)\ \to\
\R/\Z$$
is a quasi-null function.
\endproclaim

This is clear since Theorem 7.6 says
$$\overline\rho(M)\cdot\alpha-\overline\rho(M)\cdot\beta\ =\
\tilde\rho(M)\cdot\alpha-\tilde\rho(M)\cdot\beta$$
if $\alpha$ and $\beta$ can be connected by a path in
$\RR_k(\pi)$
and
$\overline\rho(M)=\tilde\rho(M)=0$ on the trivial
representation.

\proclaim{10.5. Corollary} If $M$ and $M^\prime$ are oriented
homotopy
equivalent, then $\overline\rho(M)\ =\
\overline\rho(M^\prime)$ up
to a
quasi-null function. Moreover, if $\alpha$ factors through a
group
satisfying the
Novikov conjecture, then $\overline\rho(M)\cdot\alpha\ -\
\overline\rho
(M^\prime)\cdot\alpha\ \in\ \Q/\Z$.
\endproclaim

Recall that the Novikov conjecture for a group $\pi$ asks that
for
any
homomorphism $\theta:\pi_1(M)\to\pi$,\ $M$ a closed oriented
manifold,
the homology class $\theta_\ast(l(M))\in H_\ast(\pi;\Q)$
depends
only on
the oriented homotopy class of $M$\ ($l(M)$ is the total
homology
class
dual to the Hirzebruch class $L^\ast(M)$). Now
$\alpha\mapsto\overline\rho_
{\alpha\theta}(M)$ defines a continuous function
$\RR_k(\pi)\to \R/\Z$ \ (if $M$ is odd-dimensional). Consider
the
subspace
$\RR_k^0(\pi)\subseteq \RR_k(\pi)$, a union of some connected
components,
consisting of all $\alpha$ whose associated flat bundle
$\xi_\alpha$
over
$B\pi$ is trivial on all finite subcomplexes of $B\pi$. Now,
following
\cite{\APS}, part III, each $\alpha\in \RR_k^0(\pi)$ defines
an
element
$v_\alpha\in \ K^{-1}(B\pi)\otimes\R/K^{-1}(B\pi)$ using a
path of
connections
from $\alpha$ to a trivial connection (i.e. trivial
monodromy), and
then
we have the formula:
$$\overline\rho_{\alpha\theta}(M)\ \equiv \
<\theta^\ast\ch(v_\alpha)
\cup\L^\ast(M),[M]>\ \mod \Z$$
where $4^r\L^{4r}(M)=L^{4r}(M)$ for each nonnegative integer
$r$. If
$\pi$
satisfies the Novikov conjecture, then the right side is a
homotopy
invariant
of $(M,\theta)$. For general $\alpha$ we have non-trivial
$\xi_\alpha$,
but, since it is a flat bundle, all the Chern classes vanish
and
some
multiple $m\xi_\alpha$ is (stably) trivial on any finite
subcomplex.
Thus $m\overline\rho_{\alpha\theta}(M)\ =\
\overline\rho_{(m\alpha)\theta}(M)$
is a homotopy invariant of $(M,\theta)$ and the second
sentence of
the
corollary follows.

10.6. We now return to the unreduced $\rho$-invariant. We
recall
theorem A with
some new notation. Let $\gamma$ be a curve in $\RR_k(\pi)$,
which
is holomorphic at 0, where $\pi=\pi_1(M)$ as usual. Set
$$\lambda_M(\gamma)\ =\ \lim_{t\to 0^+}\rho(M)\cdot\gamma(t)\
-\
\rho(M)\cdot\gamma(0)$$
If $\gamma$ is induced by a holomorphic curve of flat
connections,
then
Theorem 1.5 asserts that $\lambda_M(\gamma)$ is a sum of
signatures
associated
to a linking pairing defined on the torsion submodule of the
homology of
$M$ with a twisted coefficient system of $\C[[t]]$-modules
defined
by
$\gamma$. Since this depends only on duality and cup-product
in $M$
,
theorem A gives an explicit homotopy theoretic formula for
$\lambda_M(\gamma)$. In fact a recent preprint of B.Fine,
P.Kirk and
E.Klassen \cite{\FKK} shows that any holomorphic $\gamma$ can
be
lifted
to a holomorphic curve of flat connections.

We now state the main result of this section.

\proclaim{10.7. Theorem} \roster
\item The function $\rho(M):\ \RR_k(\pi)\ \to \R$ is uniquely
determined, up to a qusi-null function with values in $\Z$,
by $\lambda_M$ and $\overline\rho(M)$. Therefore,
$\rho(M)$ is uniquely determined up to a quasi-null function
by
$\lambda_M$
and $\tilde\rho(M)$.
\item
If $M$ and $M^\prime$ are oriented and homotopy
equivalent, then the difference
$$\rho(M)-\rho(M^\prime):\ \RR_k(\pi)\ \to\ \R$$
is a qusi-null function.
If $\alpha$ factors through a group satisfying the Novikov
conjecture, then
$\rho(M)-\rho(M^\prime)$ admits a rational value on the
component of
the
representation space $\RR_k(\pi)$ containing $\alpha$.
\endroster
\endproclaim

Theorem 10.7 will follow from the following elementary
consequence
of the
"curve selection lemma" (see \cite{\BC}, prop. 8.1.17):
\proclaim {10.8. Lemma} Let $V$ be a real algebraic set, and
$\phi:V\to\R$
a function satisfying:
\roster
\item $\phi $ is piecewise-continuous, i.e. there is a
stratification of
$V$ by subvarieties $V=V_0\supseteq V_1\supseteq V_2\dots
\supseteq
V_n=
\varnothing$ such that $\phi|V_i-V_{i+1}$ is continuous;
\item the reduction $\overline\phi: V\to \R/\Z$ is locally
constant;
\item for any curve $\gamma$ in $V$ which is holomorphic at 0,
$$\lim_{t\to 0^+}\phi(\gamma(t))\quad=\quad\phi(\gamma(0)).$$
\endroster
Then $\phi$ is locally constant.
\endproclaim
\demo{Proof} We may assume, by (2), that $\overline\phi=0$
after
adding a locally constant function. Thus $\phi(V)\subset\Z$
and
we only need to show $\phi$ is continuous. Suppose that
$\phi|V_{i+1}$
is continuous and $x\in V_i$ is a discontinuity of $\phi|V_i$.
By
(1) we
conclude that $x\in V_{i+1}$. Let $\{x_n\}$ be a sequence of
points
in
$V_i$ converging to $x$ such that $\phi(x_n)\ne\phi(x)$ for
all $n$.
Since $\phi|V_{i+1}$ is continuous, we may assume
$\{x_n\}\subseteq
V_i-V_
{i+1}$. Now $V_i-V_{i+1}$ has a finite number of components,
each of
which
is a semi-algebraic set (\cite{\BC}, th. 2.4.5) and so we may
assume
that
$\{x_n\}$ is contained in one of them $C$. Since $\phi|C$ is
continuous,
it is constant, and so $\phi(C)\ne\phi(x)$. Now we apply the
curve
selection
lemma to obtain a holomorphic curve $\gamma$ in $V_i$ such
that
$\gamma(0)=
x$ and $\gamma(t)\in C$ for all $t\in (0,\epsilon)$ for
some $\epsilon >0$. But now
(3) contradicts $\phi(C)\ne\phi(x)$.
\enddemo

One can conjecture that the quasi-null functions of Theorems
10.4
and
10.7(1) are $\Q$-valued.

In the Appendix, written by S.Weinberger, it is proved that
the quasi-null function of Corollary 10.5 and of
Theorem 10.7.(2) are $\Q$-valued.

\vfil\break

\hsize=5in
\vsize=7.8in
\strut
\vskip5truemm
\font\small=cmr8
\font\itsmall=cmti8

\def\smallarea#1{\par\begingroup\baselineskip=10pt#1\endgroup\
par}
\def\abstract#1{\begingroup\leftskip=5mm \rightskip 5mm
\baselineskip=10pt\par\small#1\par\endgroup}

\def\refer#1#2{\par\begingroup\baselineskip=11pt\noindent\left
skip=8
.25mm\rightskip=0mm
 \strut\llap{#1\kern 1em}{#2\hfill}\par\endgroup}

\nopagenumbers
\centerline{\it Appendix}
\vskip 25.3pt
\centerline{\bf RATIONALITY OF $\rho$-INVARIANTS}
\vskip 25.3pt

\centerline{SHMUEL WEINBERGER
\footnotemark"$^1$"}
\footnotetext"$^{1}$"{The research was supported by the NSF}
\vskip 15 pt
\smallarea{%
\centerline {\itsmall Department of Mathematics,}
\centerline {\itsmall The University of Chicago}
\centerline {\itsmall Chicago, Illinois}}
\vskip 0.5in
fibred
knot,



\redefine\F{{\frak F}}

\define\BCR{1}
\define\CG{2}
\define\FL{3}
\define\FW{4}
\define\KS{5}
\define\W{6}
\topmatter
\endtopmatter
\TagsOnRight
\document

Here we observe some rationality statements that follow from
\cite{\FL} and known facts regarding the Novikov conjecture.
\proclaim{Theorem}
\roster
\item"{(i)}" If $h:M^\prime \to M$ is a homotopy equivalence,
then $\rho(M^\prime)-\rho(M)$ take values in $\Q$.
\item"{(ii)}" If a nontrivial finite group acts freely and
homologically trivially on $M$, then $\rho(M)$ takes values in
$\Q$.
\endroster
\endproclaim

In (ii) we assume that that a nontrivial finite group $G$ acts
freely on $M$
such that the sequence
$1\to \pi_1(M)\to \pi_1(M/G)\to G\to 1$ splits and the action
of $G$
on
$H_\ast(M;\Q[\pi_1(M)])$ is trivial.

Case (ii) is an analogue of a conjecture
of Cheeger and Gromov \cite{\CG} in light of the existence of
"F-
structure"
choppings of complete manifolds with bounded curvature and
finite
volume.

\demo{Proof} (i) follows from \cite{\FL}, Corollary 10.5 and
the
following
propositions:
\enddemo
\proclaim{Proposition 1} The Novikov conjecture is correct for
$\Gamma\subset GL_n(\overline \Q)$ (where $\overline \Q$ is
the
algebraic
closure of $\Q$).
\endproclaim
This is proven in \cite{\KS} and \cite{\FW}.
\proclaim{Proposition 2} If $\pi$ is a finitely presented
group then
every
component of $\RR_k(\pi)$ contains a point defined over
$\overline\Q$. The
corresponding representation is a homomorphism
$\rho:\pi\to U_k(\overline\Q)$.
\endproclaim
\demo{Proof} Indeed, $\RR_k(\pi)$ is a real algebraic variety,
defined
over $\Q$. It is the case that the $\overline\Q$-points of any
real
variety defined over $\Q$ are dense in the $\R$-points! (Since
the
variety
is triangulable, density implies that there are points in each
component).
According to the Tarski - Seidenberg Theorem (see e.g.
\cite{\BCR})
a system
of $\Q$-equations and inequalities has an $\R$-solution iff it
has a
$\overline\Q$-solution. If $p=(p_1,\dots,p_n)$ is an $\R$-
point in
coordinates, let $r_i,q_i$ be any rational numbers with
$r_i<p_i<q_i$. There is a $\R$-point satisfying  the equations
of
the
variety and the inequalities $r_i<p_i<q_i$ (namely $p$) so
there is
such a $\overline\Q$-point, which is exactly density.
$\square$
\enddemo

To prove (ii) one relies on \cite{\W} which shows that the
classes
in oriented
bordism $\Omega(B\pi)\otimes\Q$ represented by manifolds with
homologically
trivial action are the same as those represented by
differences of
homotopy
equivalent manifolds, and (the obvous fact) that
$\rho\mod \Q$ only depends on
the class in $\Omega(B\pi)\otimes\Q$.
$\square$

\end